\documentclass[11pt,english]{article}
\usepackage{fix-cm}
\usepackage{xr}
\usepackage{makecell}
\usepackage{tabularx}
\usepackage{subcaption}
\usepackage{caption}
\usepackage{graphicx}
\usepackage{threeparttable}
\usepackage{ulem} 
\usepackage{xcolor}
\usepackage{setspace}
\usepackage{natbib}
\usepackage{booktabs}
\usepackage{ulem}
\usepackage{makecell,rotating,multirow,diagbox} 
\usepackage{amssymb}
\usepackage[mathscr]{euscript}
\usepackage{mathrsfs}
\usepackage{amsmath,bm}
\newtheorem{result}{Result}
\usepackage{amsbsy}
\usepackage{setspace}
\usepackage{soul}
\usepackage{array}
\usepackage{booktabs}
\usepackage{enumerate}
\usepackage{enumitem}

\usepackage{pifont}
\newcommand{\cmark}{\ding{51}}
\newcommand{\xmark}{\ding{55}}
\usepackage{longtable}
\usepackage{setspace}
\usepackage{appendix}
\usepackage{geometry}
\usepackage{authblk}
\usepackage{url}
\emergencystretch=\maxdimen
\usepackage[english]{babel}
\usepackage[font=normalsize, labelfont=bf]{caption}
\usepackage[font=normalsize, labelfont=bf]{subcaption}
\usepackage{natbib}
 \bibpunct[, ]{(}{)}{,}{a}{}{,}
 \def\bibfont{\small}

\usepackage{float}
\begin{document}

\title{Optimizing Urban Electric Vehicle Charging and Battery Swapping Infrastructure: A Location-Inventory-Grid Model}

\author[a]{Wenqing Ai\thanks{These authors contributed equally to this work.}}
\author[b]{Hanyu Cheng\protect\footnotemark[1]}
\author[b]{Wei Qi\thanks{Corresponding author: qiw@tsinghua.edu.cn}}
\affil[a]{School of Economics \& Management, University of Science and Technology Beijing, Beijing, 100083, China}
\affil[b]{Department of Industrial Engineering, Tsinghua University, 30 Shuangqing Road, Beijing 100084, China}
\title{Optimizing Urban Electric Vehicle Charging and Battery Swapping Infrastructure: A Location-Inventory-Grid Model}

\date{}
\maketitle

\begin{abstract}
    The rapid rise of electric vehicles (EVs) places unprecedented stress on both urban mobility systems and low-voltage power grids. Designing battery swapping and charging networks that are cost-efficient, grid-compatible, and sustainable is therefore a pressing yet complex challenge: service providers must jointly optimize station locations, battery inventory, and grid interaction under high-dimensional uncertainty. We develop an integrated location-inventory-grid model and employ a continuous approximation approach to overcome the intractability of discrete formulations. Our analysis compares centralized versus decentralized charging, with and without participation in frequency regulation. The results reveal that centralized charging, when combined with frequency regulation, not only reduces cost but also strengthens grid stability. However, it may constrain operational flexibility near the optimum, potentially challenging efforts to mitigate environmental impacts by lowering battery inventories. These results offer actionable guidance for cost-efficient, environmentally sustainable, and grid-compatible scaling of urban EV infrastructure to meet the demands of large-scale EV adoption.
    
\end{abstract}

\clearpage

\section{Introduction}

The electrification of transportation has emerged as a cornerstone of global efforts to combat climate change and modernize mobility systems. Electric vehicles (EVs) offer compelling environmental and economic advantages, such as reducing greenhouse gas emissions, improving energy efficiency, and lowering fuel costs \citep{EPA2025EVMyths}. Governments worldwide have implemented ambitious policies to promote EV adoption recognizing these benefits, with countries like the United Kingdom, Canada, and Japan announcing timelines to phase out internal combustion engine vehicles \citep{EC2021ICEPhaseOut}. As a result, global electric car sales exceeded 17 million in 2024, accounting for over 20\% of all car sales \citep{IEA2025GlobalEV}. Riding on this momentum, battery charging and swapping systems, which refuel electric vehicles by replacing depleted batteries with charged ones, are rapidly evolving. They serve as an alternative to plug-in charging due to advantages such as faster service, reduced space requirements, and enhanced safety \citep{qi2023scaling}. For example, NIO, a Shanghai-based electric vehicle startup, offers its customers lifelong battery swapping service free of charge and had deployed more than 2,600 swapping stations across China by 2024 \citep{NIO2024Video}. The eCharge4Drivers project, funded by the European Commission, has expanded its electric vehicle energy supply services across 12 countries and is also building battery swapping networks, aiming to provide seamless and efficient energy replenishment solutions throughout Europe \citep{cordis2023echarge4drivers}.

Despite the rapid growth and widespread enthusiasm surrounding battery swapping services, scaling up these networks in urban settings presents multifaceted challenges. The highly decentralized layout of swapping stations, which is necessary to ensure convenient user access without detours, significantly complicates daily operations. This dispersion amplifies unpredictable and variable local demand. As a result, service providers are forced to maintain large inventories of fully charged batteries at each station. This strategy is expensive and also raises environmental concerns due to resource consumption and battery stockpiling. At the same time, the rapid increase in localized battery charging places substantial strain on existing low-voltage power distribution grids. Many of these grids were not originally designed for such heavy loads, therefore limiting the achievable charging speeds. This infrastructural bottleneck poses a significant risk to the grid, while costly grid upgrades remain difficult to implement in densely populated areas. Moreover, although coupling EV infrastructure with power grids offers promising opportunities to enhance grid flexibility and stability through ancillary services \citep{IEA2021GridIntegration}, the precise technical and economic pathways to realize these benefits remain largely unexplored. These pathways include considerations such as site selection and inventory management of battery swapping systems. Finally, the complexity inherent in planning and optimizing large-scale swapping and charging networks such as balancing station siting, battery inventory, and grid interactions, results in computationally intensive problems. The computational problems challenge current methodologies and limit actionable managerial insights. In summary, addressing these inventory, station siting, and grid-integration challenges calls for innovative, scalable, and grid-aware solutions that integrate various aspects of infrastructure and operation to pave the way toward the EV era.

This paper aims to address these complex site-selection, inventory management, and grid-interaction challenges involved in scaling up battery swapping services in urban areas. To tackle these challenges, this paper develops a comprehensive location-inventory-grid coupling modeling framework that integrates both network configuration and grid interaction considerations. Regarding network configuration, we consider a ``swap locally, charge centrally" strategy. In this approach, batteries are exchanged at decentralized stations but transported to and charged at centralized hubs connected to high-voltage grids. This strategy offers a promising path to address these challenges, as supported by recent findings in \cite{qi2023scaling}. By pooling battery demand across multiple swapping stations and shifting the charging process to centralized facilities connected to the high-voltage grid, this strategy significantly reduces the total number of batteries required in the system while alleviating the load on distributed low-voltage power networks. This strategy was formally proposed by the State Grid Corporation of China in 2011 as ``battery swap as the mainstay, plug-in charging as the supplement, centralized charging, and unified distribution," and has been piloted in projects such as the Gaoantun station in Beijing and the Xuejiadao station in Qingdao \citep{nea2012EVstation,nea2013QingdaoEVnetwork}, and adopted in cities like Hangzhou with large centralized charging hubs \citep{gov2011EVmode,xu2015architecture}. Building on these practices, our study explores this centralized charging network's potential to alleviate inventory and charging power constraints inherent in the distributed network.

In addition to network configuration, the framework incorporates grid interaction by modeling the participation of battery swapping and charging facilities in frequency regulation services. Frequency regulation refers to a critical ancillary service that stabilizes power grid frequency by continuously adjusting electricity consumption or supply in response to real-time imbalances between demand and generation. Traditionally provided by large-scale generators, this service increasingly relies on distributed energy resources with fast-response capabilities. Battery swapping and charging stations, with their substantial and flexible battery inventories, are well-positioned to support this function. For example, in 2024, NIO Power partnered with CSG Energy Storage to integrate battery swap stations into virtual power plant platforms for frequency regulation \citep{nio2024_csg_energy_storage}. Additionally, in Denmark and the Netherlands, NIO’s Power Swap Stations serve as grid-stabilizing assets that help maintain frequency control \citep{nio2024_power_revolutionizes_ev_mobility}. Compared with other V2G or B2G schemes, such as peak shaving, frequency regulation services are more profitable \citep{white2011using}. In addition, it requires minimal battery cycling, which has a negligible impact on the state of charge and does not disrupt the primary operation of the swapping stations. By enabling frequency regulation participation, these facilities can not only enhance grid stability but also generate additional revenue streams that help offset the significant costs of infrastructure deployment and operation. Our study examines how battery swapping and charging networks can benefit from participating in frequency regulation and contribute to grid stability. To this end, our model framework includes four distinct configurations that combine decentralized and centralized network layouts with and without frequency regulation participation.

Furthermore, recognizing the computational complexity of planning large-scale swapping and charging networks that couple location decisions, inventory management, and grid interactions, we innovatively apply a continuous approximation approach. This technique transforms the original discrete, high-dimensional problem, such as the one proposed in \cite{qi2023scaling}, into a continuous cost function, allowing for efficient approximate solutions while retaining key managerial insights. By comparing these four models across multiple dimensions, our study offers practical guidance for policymakers and service providers aiming to implement scalable, grid-aware battery swapping infrastructures. 

We summarize our contributions as follows:

\begin{enumerate}
    \item To the best of our knowledge, this study is the first to address the challenges of large battery inventory and grid strain in urban battery swapping with revenue opportunities from frequency regulation. Previous studies have examined centralized charging or grid participation separately, but none have captured their combined effects. However, such combined effects have strong potential to reduce the high costs of building and operating battery swapping services and to alleviate inventory burdens and grid stress. These effects also open up new revenue streams through participation in frequency regulation. We bridge this gap by developing an integrated location-inventory-grid model. This model enables a comprehensive analysis of how infrastructure design and grid interaction jointly affect costs, environmental impact, and operational feasibility.
    
    \item Methodologically, we apply the continuous approximation (CA) approach to address the computational challenges of the location-inventory-grid coupling problem. This coupling model is highly complex in the discrete site-selection context due to its spatial and combinatorial complexity. However, by reformulating the problem in a continuous setting, we avoid the complexity of discrete site selection. This substantially reduces the computational burden while preserving the essential structural trade-offs. As a result, our approach enables efficient exploration of large-scale planning scenarios that would otherwise be computationally intractable, and provides a new methodological perspective for EV infrastructure planning.
    
    \item We conduct extensive numerical experiments to evaluate four operational configurations, defined by charging structure (centralized vs. decentralized) and frequency regulation participation (with vs. without). We compare these configurations in terms of economic cost, environmental impact, and contribution to grid stability. We find that combining centralized charging with participation in frequency regulation achieves the best performance in terms of cost-effectiveness and grid stability. At the same time, this combination maintains relatively low battery inventories. This advantage persists under forward-looking scenarios with scaling up battery swapping services and anticipated technological advancements, such as higher charging power and declining battery costs. We also find that participation in frequency regulation increases the sensitivity of cost density to inventory strategy. These results provide actionable managerial insights to support the development of cost-effective, environmentally sustainable, and grid-compatible battery swapping systems.
\end{enumerate}

The remainder of this paper is organized as follows. Section~\ref{sec:literature_review} reviews the relevant literature. Section~\ref{sec:problem_setting} presents the problem description and model setup. Section~\ref{sec:numerical_analysis} reports the results of numerical experiments and discusses the derived managerial insights. Finally, Section~\ref{sec:conclusion} concludes the paper.

\section{Literature Review}
\label{sec:literature_review}

Battery swapping services have been growing as a critical component of sustainable urban transportation. In terms of swapping and charging modes, existing studies primarily focus on the Distributed Battery Swapping and Charging System (Distributed BSCS), which include research on station location and EV routing \citep{qian2024alternating,rahmanifar2024two,li2025battery}, inventory management within battery swapping stations \citep{sun2019robust,hua122023inventory,wang2024inventory}, and optimal charging strategies \citep{sun2017optimal,sun2019robust,gull2024multi}. In general, BSCS involves decentralized swapping and charging at each station, while Centralized Charging with Battery Dispatch (CCBD) decouples these processes by centralizing battery charging and dispatching fully charged batteries to distributed stations. In contrast to the extensive literature on BSCS, the CCBD mode is less studied in the literature. Among previous literature, \cite{ni2020inventory} propose a two-stage optimization framework for a geographically distributed battery-swapping station network including centralized charging and decentralized swapping, integrating long-term inventory planning and short-term vehicle routing to maximize network revenue with provable performance guarantees. \cite{he2022optimal} propose a planning framework for centralized EV battery swapping stations based on spatial-temporal load forecasting, integrating road topology, user travel patterns, PV deployment, and distribution grid constraints to minimize annual system costs. \cite{shaker2023joint} propose a programming model consisting of linear and nonlinear parts for battery swapping stations with centralized charging to optimize station placement, battery logistics, and grid reinforcement while minimizing total costs. Building on this foundation, this study explores the operational implications of CCBD in comparison with the traditional BSCS framework. Specifically, we develop integrated location-inventory models under both the CCBD and traditional distributed BSCS frameworks, aiming to identify optimal decision-making strategies for system service providers. Our findings provide new insights into cost efficiency and environmental impact, highlighting the potential advantages of the CCBD model.

In addition, existing studies have highlighted the potential of battery swapping and charging stations (BSCSs) to participate in Battery-to-Grid (B2G) interaction \citep{cui2023operation}. In general, BSCS providers can participate in B2G interactions and generate revenue through two main approaches: discharging batteries to sell electricity for energy trading, and providing frequency regulation services to earn ancillary service income. For the discharging approach, the comprehensive optimization model proposed by \cite{shaker2023joint} incorporates B2G discharging revenue as a key component of the objective function, significantly impacting the overall cost of establishing and operating battery swapping services. \cite{ding2022joint} demonstrate that integrating B2G discharge in orderly charging schedules reduces grid load fluctuations and increases BSS provider revenue, through Monte Carlo simulation and multi-objective optimization. BSCSs also have strong potential to participate in frequency regulation, with their large-scale battery inventories and flexible charging/discharging capabilities similar to conventional controllable loads like industrial refrigeration units and boilers\citep{zarkoob2013optimal}. For instance, \cite{wu2021two} propose a two-stage self-scheduling model for battery swapping stations that incorporates a data-driven battery degradation model under frequency regulation. \cite{zhang2021operation} propose a two-stage operation strategy for electric vehicle battery swap station clusters to participate in frequency regulation services while ensuring swap demand, thereby improving overall operational economy. \cite{wang2024optimal} propose a day-ahead bidding and real-time scheduling strategies for battery swapping stations to participate in frequency regulation, combining an information-gap decision theory (IGDT) based day-ahead bidding model with a real-time response strategy to maximize profits under regulation signal uncertainty. However, the integration of frequency regulation revenue with the location-inventory decision-making of BSCSs remains underexplored. Table~\ref{tab:literature_compare} summarizes existing studies on BSCS planning and its integration with frequency regulation, highlighting this research gap. This study addresses this gap by incorporating frequency regulation profit functions into the optimization models of both centralized charging with battery dispatch and traditional distributed BSCS modes, aiming to reveal the influence of frequency regulation participation on location and inventory strategies, as well as on grid stability. This integrated perspective offers a new direction for mitigating the high construction costs associated with BSCSs.

\begin{table}[htbp]
\small
\centering
\caption{Research comparison in literature.}
\renewcommand{\arraystretch}{2.2}
\begin{tabularx}{\textwidth}{p{3.7cm}cccccc}
\toprule
\textbf{Paper} & \multicolumn{2}{c}{\textbf{Charging mode}} & \multicolumn{2}{c}{\textbf{Decision}} & \multirow{2}{*}{\textbf{Grid Coupling}} & \multirow{2}{*}{\textbf{Methodology}} \\
\cmidrule(lr){2-3} \cmidrule(lr){4-5} 
& On-site & Centralized & Location & Dispatch & &  \\
\midrule
\cite{qian2024alternating} & \cmark & \xmark & \cmark & \xmark & \xmark & ILP \& DP \\
\cite{rahmanifar2024two} & \cmark & \xmark & \cmark & \xmark & \xmark & \parbox{3cm}{\centering Multi-Objective MILP} \\
\cite{li2025battery} & \cmark & \xmark & \cmark & \xmark & \xmark & \parbox{3cm}{\centering IP} \\
\cite{sun2019robust} & \cmark & \xmark & \xmark & \xmark & \xmark & \parbox{3cm}{\centering Robust IP} \\
\cite{gull2024multi} & \cmark & \xmark & \xmark & \xmark & \xmark & \parbox{3cm}{\centering Multi-Objective MIP} \\
\cite{wang2024inventory} & \cmark & \xmark & \xmark & \xmark & \xmark & \parbox{3cm}{\centering Robust MILP} \\
\cite{hua122023inventory} & \cmark & \xmark & \xmark & \cmark & \xmark & \parbox{3cm}{\centering MILP} \\
\cite{ni2020inventory} & \cmark & \cmark & \cmark & \cmark & \xmark & \parbox{3cm}{\centering Two-Stage BILP} \\
\cite{he2022optimal} & \cmark & \cmark & \cmark & \cmark & Regulation & \parbox{3cm}{\centering SOCP} \\
\cite{shaker2023joint} & \xmark & \cmark & \cmark & \cmark & Discharge & \parbox{3cm}{\centering BILP \& NLP} \\
\cite{ding2022joint} & \xmark & \cmark & \xmark & \xmark & Discharge & \parbox{3cm}{\centering Multi-Objective MIP} \\
\cite{wu2021two} & \cmark & \xmark & \xmark & \xmark & Regulation & \parbox{3cm}{\centering Robust Two-Stage MIP} \\
\cite{zhang2021operation} & \cmark & \xmark & \xmark & \xmark & Regulation & \parbox{3cm}{\centering Two-Stage MILP} \\
\cite{wang2024optimal} & \cmark & \xmark & \xmark & \xmark & Regulation & \parbox{3cm}{\centering IGDT model} \\
\cite{qi2023scaling} & \cmark & \cmark & \cmark & \cmark & \xmark & \parbox{3cm}{\centering MIP} \\
This paper & \cmark & \cmark & \cmark & \cmark & Regulation & CA \\
\bottomrule
\label{tab:literature_compare}
\end{tabularx}
\end{table}

The deployment of battery swapping and charging stations essentially involves a joint location-inventory optimization problem, which is a classical challenge in various logistics settings. Foundational studies have addressed similar structures, such as the inventory-location model \citep{daskin2002inventory} and warehouse-retailer network design \citep{teo2004warehouse}. Previous research on such problems has focused on developing exact algorithms to solve the underlying models, including approaches like branch-and-price \citep{teo2004warehouse} and Lagrangian relaxation \citep{daskin2002inventory}, which provide theoretical optimality but often suffer from poor scalability in large-scale networks. In the context of battery swapping and charging infrastructure, \cite{qi2023scaling} develop an integrated location-inventory model for centralized battery swapping systems, introducing customized solution techniques to improve scalability. However, the computational efficiency remains constrained by the discrete nature of the modeling framework. To overcome this, continuous approximation (CA) approach can be adopted, which transforms discrete facility location and inventory problems into continuous spatial optimization over density functions, significantly reducing computational complexity in large-scale planning. This method has since been widely adopted in logistics literature and proves particularly effective for infrastructure planning. \cite{lim2017agility} apply CA methods to analyze the trade-off between proximity and agility in supply chain distribution network design under demand uncertainty. \cite{carlsson2018coordinated} apply CA methods in delivery system to evaluate the efficiency of coordinated logistics systems involving both trucks and a drone. \cite{belavina2021grocery} apply CA methods to study the impact of grocery-store density on the food waste generated at stores and by households. Inspired by the inherent simplicity and widespread application of CA methods, this study also adopts a CA-based formulation proposed by \cite{Daganzo2005}, which models the system cost of battery swapping stations as a continuous quasi-convex function, which facilitates the identification of globally optimal solutions in a straightforward manner.

\section{Problem Description}
\label{sec:problem_setting}

\subsection{Spatial Structure}

Consider battery swapping service networks spanning a continuous urban area $\Omega \subset \mathbb{R}^2$, where two types of facilities support battery services: (i) battery swapping stations that directly serve EV users' swapping demand, and (ii) battery charging stations that recharge batteries. The network operates under one of two paradigms, the decentralized charging paradigm and the centralized charging paradigm. The two paradigms differ significantly in their spatial structure. In the decentralized charging paradigm, each battery swapping station is equipped with its own charging infrastructure and independently manages battery charging and inventory. Only the swapping stations participate in frequency regulation. In contrast, in the centralized charging paradigm, the network is organized hierarchically. The centralized charging paradigm extends the decentralized layout by introducing a separate charging station that supplies multiple swapping stations. These swapping stations are spread across the area and focus solely on providing exchange services. Both centralized charging stations and swapping stations can participate in grid frequency regulation in the centralized paradigm.

\subsection{Sequence of Events}

The battery swapping service process is the same for users but differs for service providers under the two paradigms. EV drivers simply exchange depleted batteries for fully charged ones. For service providers, the decentralized paradigm requires each swapping station to manage its own charging and inventory processes. The centralized paradigm relies on both charging stations and swapping stations, with separate swapping and charging processes.

In the decentralized charging paradigm, each swapping station is equipped with its own charging infrastructure, and performs the charging of depleted batteries and frequency regulation on-site. Once a depleted battery is collected from a user, the station immediately initiates the charging process using its local equipment. To ensure service reliability, each station maintains a primary inventory level $r^B$, which serves as the minimum number of fully charged batteries that must be available at any time. The locally managed battery stock is leveraged to provide frequency regulation services by adjusting charging and discharging activity in response to grid signals.

\begin{figure}[htbp]
	\centering
	{\includegraphics[scale=0.55]{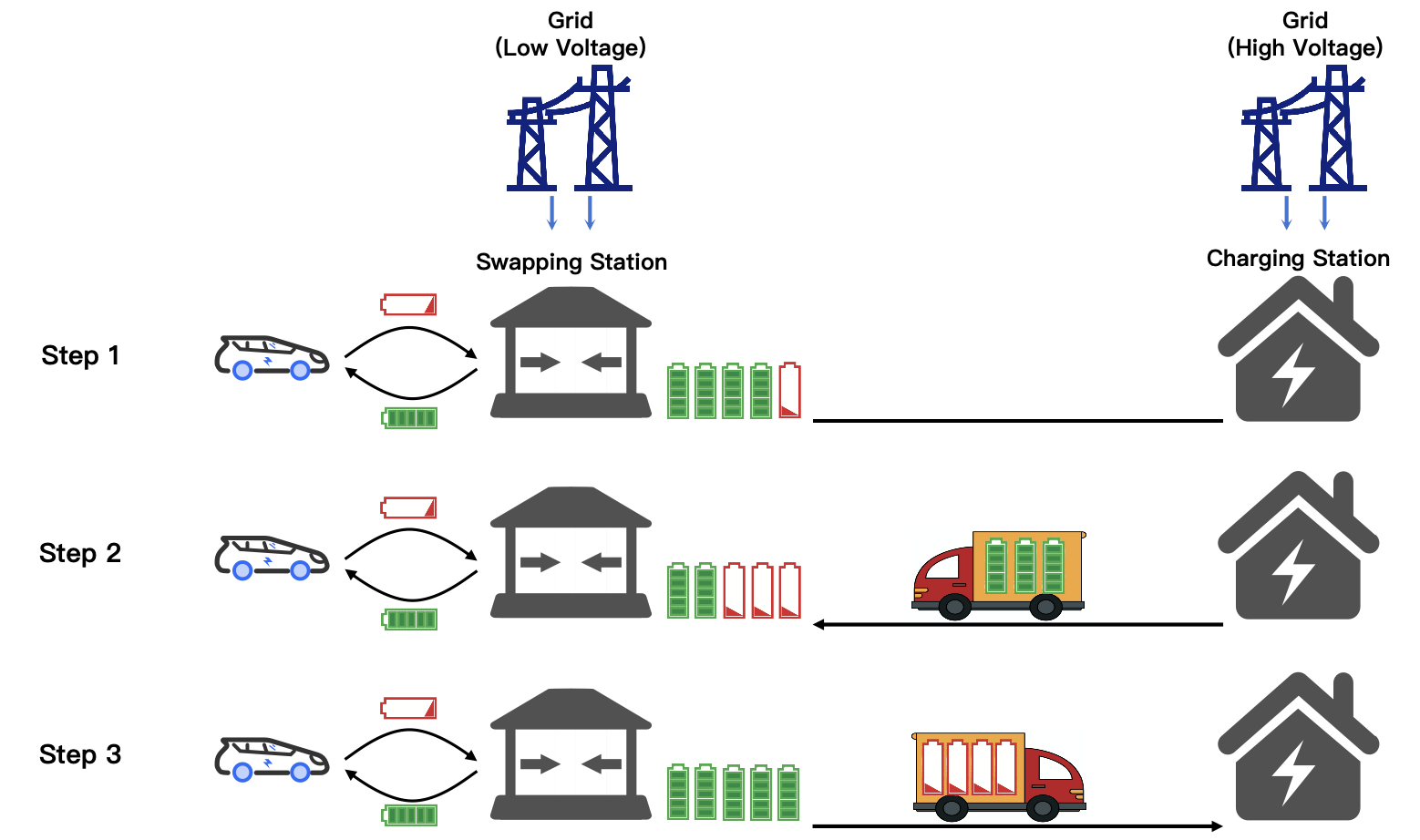}}
	\caption{
    Process of Operations in Centralized Charging Paradigm (Example of one charging station and one swapping station). \label{fig:process_of_operations}}
\end{figure}

The centralized charging paradigm builds upon this structure by introducing centralized charging stations. While swapping stations still retain their local charging capabilities, the battery charging is outsourced to centralized charging stations. Therefore, the charging operation is separated from battery swapping. The process of the operations can be seen in Figure~\ref{fig:process_of_operations}. After depleted batteries are collected at a swapping station, they are temporarily stored until the inventory level of fully charged batteries drops to a re-order point $r$ (step 1). At that moment, a transportation vehicle (e.g., a truck) is dispatched to deliver a batch of $Q$ fully charged batteries from the centralized charging station to the swapping station (step 2), and to pick up depleted batteries for return (step 3). These depleted batteries are then recharged and later redistributed to this or other swapping stations. The lead time for replenishing batteries at the swapping station equals the one-way transportation time, denoted by $T^T$. The charging station maintains a primary inventory level $R$ of fully charged batteries to ensure timely dispatch and to control the risk of stockouts at the swapping stations. This battery circulation strategy enables the service provider to maintain reliable service for EV users by mitigating inventory shortages. Moreover, with greater charging power than decentralized systems, the centralized stock of batteries offers greater capacity for participating in frequency regulation.

\subsection{Objective and Decisions}

The objective of the battery swapping service provider is to minimize the total system cost in the context of scaling up battery swapping services. For any given scale of the battery swapping station network, the total cost includes infrastructure setup costs (such as the construction of battery charging stations), battery depreciation costs, electricity costs for charging, and transportation costs incurred from relocating batteries between stations under the centralized charging paradigm. These costs are partially offset by revenues generated through participation in grid frequency regulation markets, which are enabled by the flexibility and availability of battery stock and charging infrastructure.

To achieve this objective, the service provider must make decisions regarding the locations of stations and the inventory policies for batteries. Our model builds on an existing decentralized battery swapping network and investigates how adding centralized charging stations and participating in frequency regulation affect the system. In this context, location decisions only concern charging stations, while the locations of swapping stations are considered fixed. In the decentralized charging paradigm, no location decision is required, as the swapping station network is considered already established. Inventory management focuses on maintaining a local threshold $r^B$ to ensure sufficient availability of fully charged batteries and avoid stockouts. In the centralized charging paradigm, the location decision involves selecting the placement of centralized charging stations to support the existing network of swapping stations. Inventory policies include setting the re-order point $r$ and replenishment quantity $Q$ at each swapping station, along with the primary inventory level $R$ at each charging station to ensure timely supply of recharged batteries.

\section{Optimal Stocks of Charging and Swapping Stations}

In this section, we develop continuous approximation models to determine the optimal battery stock levels in both centralized and decentralized battery swapping systems. 
Specifically, in Section~\ref{sec:cont_approx}, we introduce the continuous approximation framework that models spatially distributed swapping and charging stations, allowing us to analyze a representative local district under homogeneous conditions. 
In Section~\ref{sec:charging_ops}, we analyze the operations of charging stations and derive the aggregate primary stock required to satisfy a given service level, accounting for transportation and charging lead times as well as demand variability. 
Finally, in Section~\ref{sec:swapping_ops}, we determine the aggregate re-order points for swapping stations under an $(r, Q)$ policy for centralized charging, and the spare battery requirements for decentralized charging. 
These results provide tractable expressions for inventory management that will be used in subsequent cost and battery density analyses. The notations employed in the analysis are listed in Table~\ref{tab:notations}.

\begin{longtable}{p{2.4cm}p{13.6cm}}
\caption{Summary of Notation.} \label{tab:notations} \\
\toprule
\textbf{Notation} & \textbf{Description} \\
\midrule
\endfirsthead

\multicolumn{2}{l}{\textbf{Table \thetable{} continued from previous page}} \\
\toprule
\textbf{Notation} & \textbf{Description} \\
\midrule
\endhead

\midrule \multicolumn{2}{r}{\textit{Continued on next page}} \\
\endfoot

\bottomrule
\endlastfoot

\textbf{Sets} & \\
$\mathcal{Z}$ & Set of bidding periods, $\mathcal Z=\{1,...,Z\}$ \\
$\mathcal{W}$ & Set of time intervals, $\mathcal W=\{1,...,W\}$ \\
$\mathcal{N}$ & Set of scenarios, including peak and off-peak (idle), $\mathcal N=\{1,...,N\}$ \\

\textbf{Parameters} & \\
$c^C$ & Depreciation cost of charging station (USD \$)\\
$c^S$ & Depreciation cost of swapping station (USD \$)\\
$c^T$ & One-way transportation cost (USD \$)\\
$q$ & Truck speed ($\text{km}\cdot \text{h}^{-1}$) \\
$T^T$ & One-way truck trip time ($\text{h}$)\\
$T^C$ & The time to fully charge a depleted battery ($\text{h}$)\\
$\Delta$ & The total time for round-trip transport and charging, i.e., $\Delta = 2T^T + T^C$ \\
$B^C$ & Capacity of a battery in the system ($\text{kW}\cdot \text{h}$) \\
$c^B$ & Depreciation cost of a battery not participating in frequency regulation (USD \$) \\
$c^R$ & Depreciation cost of a battery participating in frequency regulation (USD \$) \\
$c^E_n$ & Electricity price of scenario $n\in\mathcal{N}$ ($\text{USD \$}\cdot \text{kW}^{-1}\cdot \text{h}^{-1}$) \\
$\lambda^C$ & Charging/discharging power in charging stations (kW) \\
$\lambda^S$ & Charging/discharging power in swapping stations (kW) \\
$\mu(l)$ & The average number of swap demands at a single swapping station per hour around location $l$  \\
$\sigma^2(l)$ & The variance of swap demands at a single swapping station per hour around location $l$  \\
$\bar{\mu}_{zn}(l)$ & The average number of swap demands per unit area around location $l$ in $z$th period under $n$th scenario \\
$\bar{\sigma}^2_{zn}(l)$ & The variance of swap demands per unit area around location $l$ in $z$th period under $n$th scenario \\
$\bar{\mu}_z(l)$ & The average number of swap demands per hour per unit area around location $l$ in period $z$ \\
$\bar{\sigma}_z^2(l)$ & The variance of swap demands per hour per unit area around location $l$ in period $z$ \\
$\bar{\mu}(l)$ & The average number of swap demands per hour per unit area around location $l$ \\
$\bar{\sigma}^2(l)$ & The variance of swap demands per hour per unit area around location $l$ \\
$\epsilon^S$ & The pre-specified stockout probability at swapping station \\
$\epsilon^C$ & The pre-specified stockout probability at charging station \\
$\epsilon^{BS}$ & The pre-specified stockout probability at swapping station in benchmark \\
$\epsilon^B$ & Pre-specified probability of capacity shortage \\
$\theta$ & Prescribed level of the regulation performance \\
$\eta_z$ & Minimum (stochastic) regulation capacity proportion required to fulfill the mileage requirement \\
$\kappa_n$ & Proportion of time in scenario $n \in \mathcal{N}$ \\
$\bar{p}_z$ & Mean of market clearing price for regulation capacity in period $z$ ($\text{USD \$}\cdot\text{MW}^{-1}$) \\
$\phi(Q)$ & Piecewise-defined variance function for re-order quantity \\
$\nu$ & Threshold separating two regimes in piecewise variance function $\varphi(Q)$ \\
$\Phi^{-1}(\cdot)$ & Inverse CDF of the standard normal distribution \\
$g_{w,z}$ & Random AGC signal, $g_{wz}\in[-1,1]$ \\
$\eta_z$ & Random variable whose distribution depends on the AGC signal \\
$p_z$ & Regulation market clearing price for capacity in $z$th period, which mean is $\bar p_z$ ($\text{USD \$}\cdot\text{MW}^{-1}$) \\
$L_z(l)$ & Expected number of batteries in transit in period $z$ at location $l$ \\
$\mathbb{E}[\cdot]$ & Expectation operator \\

\textbf{Variables} & \\
$R(l)$ & Amount of charged batteries that charging stations initially stocks per unit area around location $l$ \\
$r(l)$ & Aggregate re-order point of all swapping stations per unit area around location $l$  \\
$Q(l)$ & Re-order battery quantity from the charging station to the swapping station for stations around location $l$ \\
$\bar{B}_z(l)$ & Regulation capacity in $z$th period per unit area around location $l$  \\
$r^B(l)$ & Aggregate battery stock level for all swapping stations with decentralized charging per unit area around location $l$  \\
$\rho_s(l)$ & Density of swapping stations around any location $l$ within the city region $\Omega$ \\
$\rho_c(l)$ & Density of charging stations to deploy around $l$ \\

\textbf{Functions} & \\
$B_z(t,l)$ & Service provider's available capacity per unit area around location $l$  at time $t$, where $t$ belongs to $z$th period \\
$L_z(t)$ & The number of batteries being transported from swapping stations to charging stations (or from charging stations to swapping stations) per unit area around location $l$ at time $t$ in period $z$ \\
$G_z(t)$ & The number of batteries at charging stations per unit area around location $l$ at time $t$ in period $z$ \\
$S_z(t)$ & The number of batteries at swapping stations per unit area around location $l$ at time $t$ in period $z$ \\
$C(\rho_c(l), Q(l))$ & Cost density function of centralized charging around location $l$ \\
$J^{CA}$ & Total cost of centralized charging in a region \\
$B(l)$ & Battery density function of centralized charging around location $l$ \\
$\bar{B}_z(l)$ & Regulation capacity of centralized charging per unit area around location $l$ in period $z$ \\
$C^B(l)$ & Cost density function of decentralized charging around location $l$ \\
$B^B(l)$ & Battery density function of decentralized charging around location $l$ \\
$\bar{B}^B_z(l)$ & Regulation capacity of decentralized charging per unit area around location $l$ in period $z$ \\
\end{longtable}

\subsection{Continuous Approximation Method}
\label{sec:cont_approx}

Inspired by the literature on continuous approximation models for logistics systems \citep{Daganzo2005}, we model the densities of swapping stations and demand as spatially continuous functions that vary gradually across the region. Specifically, let $l \in \Omega$ denote a location in the service area. Parameters and variables such as the swapping station density $\rho_s(l)$, the mean and variance of battery swapping demand per unit area $\bar{\mu}(l)$ and $\bar{\sigma}^2(l)$, and the charging station density $\rho_c(l)$ may vary significantly across the entire urban region. Within a local service district such as a subregion covered by a single charging station, however, these quantities are approximately constant. Therefore, we drop the location index $l$ and analyze a representative local district under a homogeneous setting.

In each district, (1) each swapping station faces mean demand $\mu = \bar{\mu} / \rho_s$ and variance $\sigma^2 = \bar{\sigma}^2 / \rho_s$, (2) all swapping stations adopt a common reorder quantity $Q$, and (3) one charging station serves $\rho_s / \rho_c$ swapping stations. We adopt the Manhattan (or $L_1$) distance metric and tile the city with diamond-shaped (rotated-square) districts.  
For a district of area, the average one-way truck distance from a swapping station to the district center, i.e., the charging station, is $\int_{r=0}^{(1/2)\sqrt{1/\rho_c}} \sqrt{2} r \cdot 8r ~\mathrm{d}r \cdot \rho_c = \frac{\sqrt{2}}{3\sqrt{\rho_c}}$. Hence, the corresponding one-way travel time is $T^{T} = \frac{\sqrt{2}}{3q\sqrt{\rho_c}}$, and the transportation cost is $c^T \cdot \frac{\sqrt{2}}{3\sqrt{\rho_c}}$, where $q$ is the speed of the truck and $c^T$ is the unit cost per kilometer of transport.

\subsection{Operations of Charging Stations} 
\label{sec:charging_ops}\leavevmode

Charging stations constitute the backbone of the centralized battery management system, where depleted batteries from decentralized swapping stations are transported, recharged, and redistributed. To obtain a tractable expression for the required aggregate primary battery stock $R$ at the charging stations within a unit area, we first define $D(t)$ as the \textit{battery deficit} at a charging station at time $t$, reflecting the net shortage from all swapping stations it serves. This represents the difference between the total number of depleted batteries received and the number of fully charged batteries returned. To ensure that the stockout probability at the charging station remains below a threshold $\epsilon^C$, we impose the constraint $\text{Prob}\left(\rho_c D(t) > R \right) < \epsilon^C$. The term $\rho_c D(t)$ denotes the total deficit aggregated over all charging stations per unit area.

To quantify $D(t)$, we draw on the discrete analysis of \cite{qi2023scaling} and extend their results to a continuous model. We first derive the expected value of the battery deficit at a charging station, denoted by $D(t)$. Following Result 3 in \cite{qi2023scaling}, the expected deficit at each swapping station $i$ is given by $\mathbb{E}\left[D_i(t) \right] = (T^T+T^C)\mu + Q - m$, where $T^C$ is the time required to fully charge a battery, and $m$ can be assumed to be zero without loss of generality, according to their Result 6. Therefore, the expected deficit of a charging station is given by $\mathbb{E}[D(t)] = \left((T^T + T^C)\mu + Q\right)\frac{\rho_s}{\rho_c}$, where $\frac{\rho_s}{\rho_c}$ represents the number of swapping stations served by one charging station.

Next, we characterize the variance of $D(t)$.
Because the exact distribution of the battery deficit $D(t)$ is intractable, we adopt the piece-wise approximation developed in Result 4 of \cite{qi2023scaling}. For a single swapping station $i$, the variance of the deficit is
\begin{equation}
\phi_i(Q) = 
\begin{cases} 
\Delta \sigma^2 + \frac{Q^2 - 1}{6}, & \text{if } 0 \leq Q \leq \nu, \\
Q \Delta \mu - (\Delta \mu)^2, & \text{if } Q \geq \nu,
\end{cases}
\label{eq:phi_Q_dis}
\end{equation}
where $\Delta = 2T^T + T^C$, and $\nu$ is the threshold at which the two expressions intersect. To translate \eqref{eq:phi_Q_dis} into the continuous density framework, we replace $\sigma^2$ with $\bar\sigma^2/\rho_s$, $\mu$ with $\bar\mu/\rho_s$, and $T^T$ with $\frac{\sqrt{2}}{3q\sqrt{\rho_c}}$, and multiply each term by $\frac{\rho_s}{\rho_c}$, which corresponds to the number of swapping stations served by one charging station. We then obtain the continuous density form of $\phi(Q)$ as follows:
\begin{equation}
{\phi}(Q) = 
\begin{cases} 
\frac{\bar{\sigma}^2}{\rho_c} \left( \frac{2\sqrt{2}}{3q\sqrt{\rho_c}} + T^C \right) + \frac{Q^2 - 1}{6} \frac{\rho_s}{\rho_c}, & \text{if } 0 \leq Q \leq \nu, \\
\left( \frac{2\sqrt{2}}{3q\sqrt{\rho_c}} + T^C \right) \frac{\bar{\mu}}{\rho_c}Q - \left( \frac{2\sqrt{2}}{3q\sqrt{\rho_c}} + T^C \right)^2 \frac{\bar{\mu}^2}{\rho_s\rho_c}, & \text{if } Q \geq \nu.
\end{cases}
\label{eq:phi_Q}
\end{equation}

Invoking the central limit theorem, we approximate the aggregate deficit across all swapping stations as normally distributed. Consequently, the aggregate required primary battery stock of charging stations per unit area that guarantees a stock-out probability below $\epsilon^C$ is given by the following result:
\begin{result}
\label{result:centralized_charging_R}
In centralized charging paradigm, the aggregate required primary stock of charging stations per unit area is given by:
\begin{align}
R \geq \left( \frac{\sqrt{2}}{3q\sqrt{\rho_c}} + T^C \right) \bar{\mu} + Q\rho_s + \Phi^{-1}(1-\epsilon^C)\sqrt{{\phi}(Q)} \rho_c,
\label{ieq:R_continuous}
\end{align}
where $\Phi^{-1}(\cdot)$ is the inverse cumulative distribution function of the standard normal distribution.
\end{result}
Result \ref{result:centralized_charging_R} provides the aggregate primary stock level $R$ at charging stations per unit area. This incorporates the transportation and charging lead times, the re-order quantity $Q\rho_s$ for swap station replenishments, and an additional safety term that ensures the total battery deficit stays within the tolerance level $\epsilon^C$. The right-hand side of Inequality~\eqref{ieq:R_continuous} provides a lower bound on the primary battery stock required at each charging station to satisfy the service level constraint $\epsilon^C$. The first two terms correspond to the expected battery deficit, and the last term represents an additional stock caused by the variance of battery deficit. For the remainder of our analysis, we consider the minimum inventory level that ensures the target service level is met.

\subsection{Operations of Swapping Stations}\label{sec:swapping_ops}\leavevmode

Swapping stations serve as the front line of the battery swapping network, where electric vehicles exchange depleted batteries for fully charged ones. Efficient inventory management at these stations is crucial to ensure high service levels while minimizing the total number of batteries required. We develop the re-order $(r, Q)$ policy for swapping stations under a continuous approximation framework, where $r$ is the re-order point at swapping stations, and $Q$ is the re-order quantity. We now determine the aggregate re-order point $r$ per unit area for all swapping stations. The system starts operating at time 0 with an initial stock of $r + m\rho_s$ fully charged batteries, where $m \in [0, Q]$ represents the primary stock beyond the re-order threshold. Let $\{N^S(t), t \geq 0\}$ denote a renewal process counting the number of EVs arriving at swapping stations per unit area up to time $t$ which resets each time a re-order is triggered. To ensure that the probability of stockout during the lead time $T^T$ remains below the service threshold $\epsilon^S$, $\text{Prob}\left\{ N^S(T^T) \geq r+\rho_s\right\} \leq \epsilon^S$ must hold, where $r+\rho_s$ corresponds to the stock with one additional battery at each swapping station beyond the re-order point. As the aggregate swap demand follows a normal distribution, we have $N^S(T^T) \sim \mathcal{N}\left( T^T \bar{\mu}, \; T^T \bar{\sigma}^2 \right)$. Substituting the expression of $T^T = \frac{\sqrt{2}}{3q\sqrt{\rho_c}}$, we derive the re-order point per unit area in the following result:
\begin{result}
\label{centralized_charging_r}
In centralized charging paradigm, the aggregate re-order point of swapping stations per unit area is given by:
\begin{align}
    r \geq \frac{\sqrt{2} \, \bar{\mu}}{3q \sqrt{\rho_c}} - \rho_s + \Phi^{-1}(1 - \epsilon^S) \sqrt{ \frac{\sqrt{2} \rho_s}{3q \sqrt{\rho_c}} } \cdot \bar{\sigma}.
    \label{ieq:r_continuous}
\end{align}
\end{result}
Result~\ref{centralized_charging_r} indicates that the aggregate re-order point $r$ of swapping stations per unit area consists of the expected swapping demand during the delivery lead time and a safety buffer to ensure a service level of $\epsilon^S$. The right-hand side of Inequality~\eqref{ieq:r_continuous} provides a lower bound on the re-order point required at each swapping station to satisfy the service level constraint $\epsilon^S$. The first two terms correspond to the expected battery deficit, and the last term represents an additional stock caused by the variance of battery deficit. This expression will be used in the cost density~\eqref{eqn:C_density} and battery density~\eqref{eqn:battery_density} calculations for centralized charging in Section~\ref{sec:cost_function}.

With the above derivations, we complete the continuous approximation of inventory control policies under the centralized charging architecture. Specifically, we derive the spatially continuous form of the required stock level at charging stations to meet a prescribed service level, and develop the continuous re-order point for swapping stations under an $(r,Q)$ policy.

To provide a baseline for comparison, we consider the decentralized charging model. In the decentralized charging mode, each swapping station is equipped with charging equipment, and swapped-out batteries begin recharging immediately on-site. The random number of swap requests during the charging period is denoted by $N^S(T^C)$. To ensure that the probability of battery stockout does not exceed a threshold $\epsilon^{BS}$, the number of spare batteries $r^B$ at all swapping stations per unit area must satisfy the following probabilistic constraint $\mathbb{P}\left( N^S(T^C) \geq r^B + \rho_s \right) \leq \epsilon^{BS}$. Therefore, the swap demand during charging time follows a normal distribution $N^S(T^C)\sim\mathcal{N}\left(T^C\bar{\mu},T^C\bar{\sigma}^2\right)$.
We obtain the following condition for the minimum aggregate spare battery requirement of all swapping stations per unit area to satisfy the service level $\epsilon^{BS}$ in the following result:
\begin{result}
\label{result:r^B_continuous}
    In the decentralized charging paradigm, the aggregate spare battery stock of swapping stations per unit area is given by:
    \begin{align*}
    r^B \geq T^C\bar{\mu}-\rho_s+\Phi^{-1}(1-\epsilon^{BS})\sqrt{T^C\rho_s}\bar{\sigma}.
    \end{align*}
\end{result}
This expression is then used in the cost density \eqref{eqn:cost_density_bench} and battery density \eqref{eqn:battery_density_bench} of decentralized charging paradigm in Section~\ref{sec:cost_function}.

So far, we have derived continuous approximations for inventory management policies under both centralized and decentralized charging paradigms in battery swapping networks. Building upon this inventory analysis, we proceed in the next section to examine how battery stocks in both systems participate in frequency regulation.

\section{Optimal Regulation Capacity}

As a resource of exceptional promise for frequency-regulation markets, battery swapping and charging systems can generate considerable ancillary income. In this section, we introduce the key characteristics of the regulation market, outline the process of providing regulation services, and describe the market requirements and the concept of regulation capacity that must be satisfied. These discussions establish the basis for formulating the cost function of the charging-swapping system when participating in frequency regulation in the subsequent analysis.

\subsection{Regulation Market}\leavevmode

 We consider a day consisting of $Z$ bidding periods, indexed by $z=1,...,Z$. The day-ahead frequency regulation market operates as follows. At the start of the planning horizon ($z=0$), market participants, such as the station agency, participate in the bidding process and submit bids to the system operator (SO) for providing regulation capacity (in MWs) for each period $z\in\mathcal Z$. Market-clearing prices for 
capacity are realized after all bids have been received. Then, during each period $z\in\mathcal Z$, participants with accepted bids respond in real time to an automatic generation control (AGC) signal announced by the SO, which may call for regulation-up  (i.e., feeding power to the grid) or regulation-down (withdrawing power) at some specified rate up to the contracted capacity (the bid). The accuracy to which power dispatches track the AGC signal must fulfill a performance requirement.
 
\subsection{Bidding Process and Automatic Generation Control}\leavevmode

\begin{figure}[htbp]
    \centering
    {\includegraphics[scale=0.4]{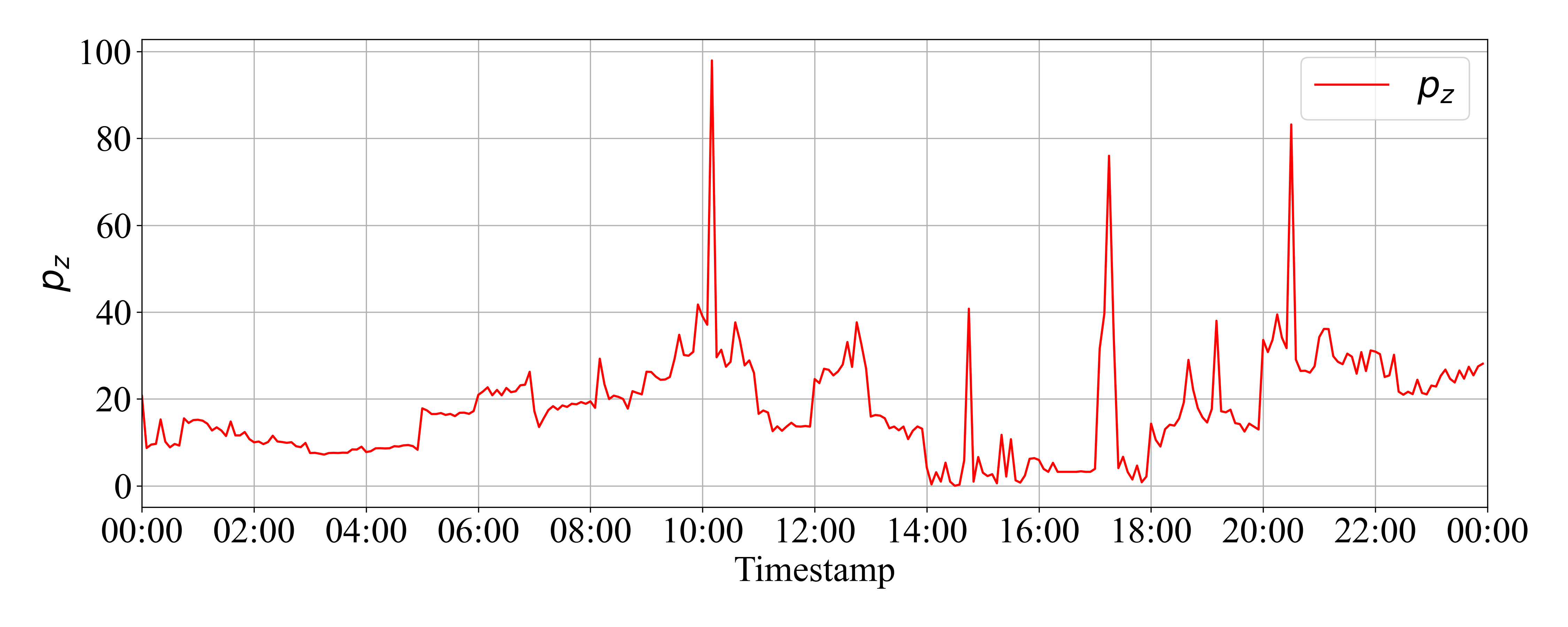}}
    \caption{Example of regulation capacity market clearing prices $p_z$ on August 31, 2024, adapted from PJM data. \label{fig:p_z}}
\end{figure}

At the beginning of the planning horizon ($z = 0$), the agency chooses to bid for a regulation capacity of $\bar B_z$ MW (i.e., an offer to provide up to $\bar B_z$ MW of contingency power to the grid on demand) for each bidding period of the day, $z \in \mathcal Z$. The market clearing price for capacity in period $z$ is denoted by $p_z$.  We consider the agency to be relatively small in capacity compared with the size of the market and thus is a price taker. Therefore, $p_z$ are modeled as random variables independent of $\bar B_z$, as their values are not known until after bids are submitted and market is cleared. An example of the variation in $p_z$ over the course of a day is shown in Figure~\ref{fig:p_z}, which is adapted from data provided by PJM\footnote{\url{https://www.pjm.com/markets-and-operations/ancillary-services}}.

\begin{figure}[htbp]
	\centering
	{\includegraphics[scale=0.4]{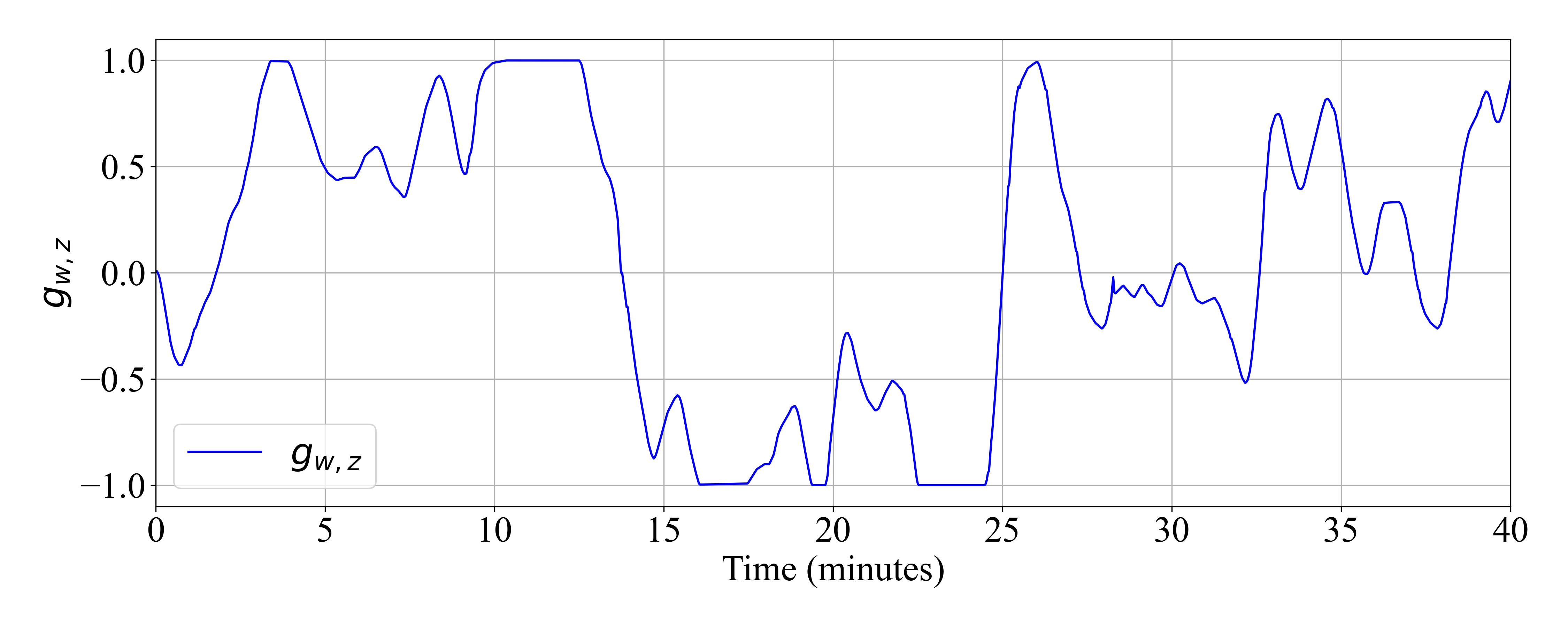}}
	\caption{
    Example of the AGC Signal, adapted from PJM data. \label{fig:agc_signal}}
\end{figure}

With a contract for capacity $\bar B_z$, the agency is obliged to dispatch power to the grid following the time-varying AGC signal during period $z$ (Figure \ref{fig:agc_signal}). In practice, the AGC signal varies 
at refined time intervals. Let each bidding period $z$ (e.g., the hour between $z-1$ and $z$ hours) be divided into a set of refined time intervals $\mathcal W$ at which the AGC signal varies (e.g., every four seconds), and let $g_{w,z}$ denote the random AGC signal. The signal is bounded within the interval [-1, 1]. During each interval $w\in\mathcal W$, the agency is programmed to automatically trace the AGC signal as closely as possible by feeding up to $\bar B_zg_{w,z}$ MW of  power to the grid (whereas a negative amount calls for regulation-down, i.e., withdrawing said power from the grid). The AGC signal is stochastic and revealed in real time, and thus are random variables at the time of bidding.

\subsection{Performance Requirement of Frequency Regulation}\leavevmode

The U.S. frequency regulation markets operate under performance-based designs.  In particular, a resource will be barred from participating in the performance-based frequency regulation market if it fails to reach a minimum performance index defined based on mileage, which measures how closely the dispatch follows the AGC signal within the bidding period. In particular, within period $z$, the total requested mileage is the cumulative change in the AGC dispatch request, i.e., $b_z=\sum_{w\in\mathcal W}|g_{w,z}\bar B_z-g_{w-1,z}\bar B_z|$. The fulfilled mileage is given by $\sum_{w\in\mathcal W}|e_{w,z}-e_{w-1,z}|$, where $e_{w,z}$ represents the resource's actual dispatch. The market design requires 
the performance index of the resource to exceed a prescribed level $\theta$, i.e., $\sum_{w\in\mathcal W}|e_{w,z}-e_{w-1,z}|\ge \theta b_z$. The performance of a resource is evaluated based on mileage (i.e., changes in dispatch) rather than actual dispatch because, to stabilize frequencies, it is more important to track the movements of the AGC signal instead of the absolute magnitude.

\begin{figure}[htbp]
	\centering
	{\includegraphics[scale=0.4]{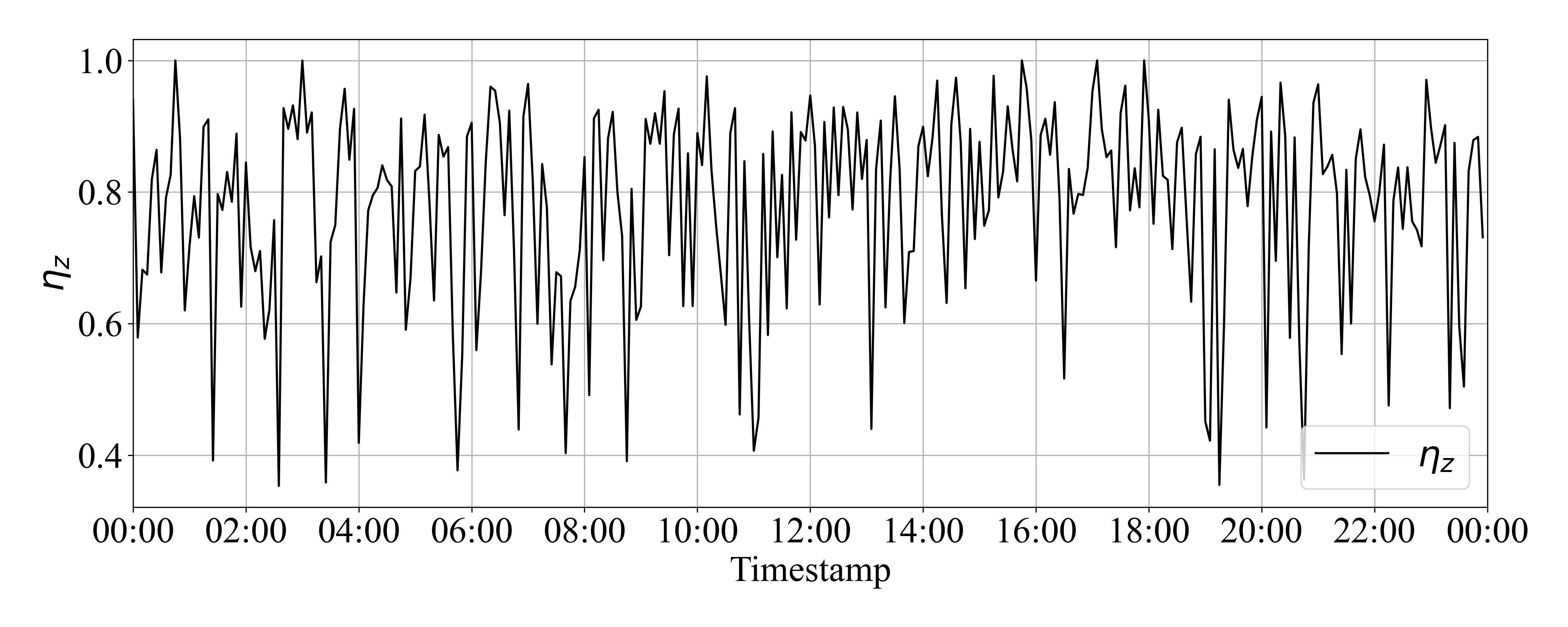}}
    \caption{Example of random variable $\eta_z$ on August 31, 2024, adapted from PJM data. \label{fig:eta_z}}
\end{figure}

Let the agency's available capacity between $w-1$th and $w$th interval in $z$th period be $B_{w,z}$, which indicates piece-wise constant available capacity in each interval (e.g., every four seconds). The maximum regulation-up (regulation-down) power that can be dispatched is capped at $B_{w,z} (-B_{w,z})$. Because dispatch is programmed to automatically follow the AGC signal that requests $g_{w,z}\bar B_z$ in $z$th interval, it is capped at $\min\{B_{w,z}, g_{w,z}\bar B_z\}$ in regulation up, and $\max\{-B_{w,z}, g_{w,z}\bar B_z\}$ in regulation down. The fulfilled mileage can be expressed as $e_{w,z}=\max\{-B_{w,z},\min\{B_{w,z},g_{w,z}\bar B_z\}\}$. Therefore, the performance requirement can be written as

\begin{equation*}
    \begin{aligned}
        \sum_{w=1}^{W}|\max\{-B_{w,z},\min\{B_{w,z},g_{w,z}\bar B_z\}\}-\max\{-B_{w,z},\min\{B_{w,z},g_{w-1,z}\bar B_z\}\}|\ge \theta b_z,
    \end{aligned}
\end{equation*}
which is not tractable in the current form due to the min and max operators associated with random variables. To rewrite this in a tractable form, we notice that the above performance requirement can be fulfilled if 
and only if the (random) available capacity is sufficiently high, given the capacity bid. Thus, we perform a variable substitution and use $\eta_{z} \bar B_z$ to denote the lowest value of $B_{w,z}$ for any given $z$ that would satisfy the performance requirement:
\begin{equation*}\begin{aligned}
\eta_{z}=&\Big\{\min \eta: \sum_{w=1}^{W}|\max\{-\eta \bar B_{z},\min\{\eta \bar B_{z},g_{w,z}\bar B_z\}\}-\max\{-\eta \bar B_{z},\min\{\eta \bar B_{z},g_{w-1,z}\bar B_z\}\}|\ge\theta b_z\Big\}\\
=&\Big\{\min \eta: \sum_{w=1}^{W}|\max\{-\eta ,\min\{\eta,g_{w,z}\}\}-\max\{-\eta,\min\{\eta,g_{w-1,z}\}\}|\ge \theta \sum_{w=1}^{W}|g_{w,z}-g_{w-1,z}|\Big\}.
\end{aligned}\end{equation*}
Note that $\eta_z$ is a random variable whose distribution depends on the AGC signal, but not the capacity bid $\bar B_z$ or other parameters. It can be interpreted as the minimum (stochastic) regulation capacity required to fulfill the mileage requirement. With this definition, the performance requirement 
can be rewritten as $B_{w,z}\ge \eta_z\bar B_z$. An illustrative example of $\eta_z$ over a day is provided in Figure~\ref{fig:eta_z}.

\subsection{Regulation Capacity}\leavevmode

To support grid frequency regulation, the service provider allocates part of the battery inventory for ancillary services. At any time $t$, the agency's available regulation capacity is determined by the total number of idle batteries at both charging and swapping stations, weighted by their respective effectiveness in providing frequency response. The available capacity in period $z$ is defined as:
\begin{equation}
    B_z(t,l)\triangleq \lambda^C G_z(t,l) + \lambda^S S_z(t,l),
    \label{eq:B_z_t_l}
\end{equation}
where $\lambda^C$ is the charging power at charging stations, $\lambda^S$ is the charging power at swapping stations, $G_z(t,l)$ and $S_z(t,l)$ represent the number of batteries stored at charging stations and at swapping stations per unit area around $l$ at time $t$ in period $z$, respectively. Due to stochastic battery demand and transport operations, $G_z$ and $S_z$ are random variables. To ensure that sufficient regulation capacity is available with high probability, the grid operator imposes a performance constraint: the probability that capacity falls below a required threshold $\eta_z \bar{B}_z$ should not exceed a tolerance level $\epsilon^B$. That is, 
\begin{equation}
    \lim_{t \to \infty} \mathbb{P}(B_z(t,l) \le \eta_z \bar{B}_z(l)) \le \epsilon^B.
    \label{ieq:B_z_B_bar_z}
\end{equation}
To obtain the steady-state distributions of $B_z(t,l)$, we use $L_z(t,l)$ to denote the number of batteries being transported from each swapping station served by the charging station to the corresponding charging station at time $t$ per unit area around $l$ in period $z$. 
The number of batteries at swapping stations around $l$ at time $t$ in period $z$ can be denoted by $S_z(t,l)=r(l)+\rho_s(l) Q(l)-L_z(t,l)$. The number of batteries at charging stations around $l$ at time $t$ in period $z$ can be denoted by $G_z(t,l)=R(l)-L_z(t,l)$. Therefore, the service provider's available capacity can be derived from Equation~\eqref{eq:B_z_t_l} and expressed as 
\begin{equation}\begin{aligned}
B_z(t,l)=\lambda^C R(l)+\lambda^S(r(l)+\rho_s Q)-(\lambda^C+\lambda^S) L_z(t,l).
\label{eq:B_z_t}
\end{aligned}\end{equation}
By definition, $L_z(t,l)$ follows a Bernoulli distribution, 
$\text{Prob}\Big(L_z(t,l)=\rho_s(l) Q(l)\Big)=\frac{T^T\bar{\mu}_{z}(l)}{\rho_s(l)Q(l)}$ and $\text{Prob}\Big(L_z(t,l)=0\Big)=1-\frac{T^T\bar{\mu}_{z}(l)}{\rho_s(l)Q(l)}$, where $\bar{\mu}_{z}(l)$ is the mean swap demand per hour per unit area around location $l$ in period $z$. Thus, the mean and variance of $L_z(t,l)$ in the steady state are $\mathbb E[L_z(l)]=T^T\bar{\mu}_{z}(l)$ and $\text{Var}(L_z(l))=T^T\bar{\mu}_{z}(l)(Q(l)-T^T\bar{\mu}_{z}(l))$. Substituting $B_z(t,l)$ from Equation~\eqref{eq:B_z_t} into Inequality~\eqref{ieq:B_z_B_bar_z}, we obtain the equivalent form
\begin{equation*}
    \lim_{t\rightarrow \infty} \Pr\Bigg( L_z(t,l) \ge \frac{\lambda^C R(l) + \lambda^S \big(r(l)+\rho_s(l)Q(l)\big) - \eta_z \bar B_z(l)}{\lambda^C + \lambda^S} \Bigg) \le \epsilon^B.
\end{equation*}
Rearranging yields the following closed-form upper bound on the available regulation capacity $\bar{B}_z(l)$ per unit area that can be offered while keeping the shortage probability below $\epsilon^B$:
\begin{result}
\label{FR_capacity}
The regulation capacity of centralized charging paradigm in $z$th period needs to satisfy 
\begin{equation}
    \begin{aligned}
        \bar B_z(l)\le \frac{\lambda^CR(l)+\lambda^S(r(l)+\rho_s(l)Q(l))-(\lambda^C+\lambda^S)\Big(\mathbb E[L_z(l)]+\Phi^{-1}(1-\epsilon^B)\sqrt{\text{Var}(L_z(l))}\Big)}{\eta_z},
        \label{ieq:reg_capacity_continuous}
\end{aligned}\end{equation}
where $\mathbb E[L_z(l)]=T^T\bar{\mu}_{z}(l)$ and $\text{Var}(L_z(l))=T^T\bar{\mu}_{z}(l)(Q(l)-T^T\bar{\mu}_{z}(l))$.
\end{result}
The numerator on the right-hand side of Inequality~\eqref{ieq:reg_capacity_continuous} captures the total available regulation power in the system. Specifically, the first term $\lambda^C R(l)$ represents the maximum power capacity contributed by centralized charging stations, and the second term $\lambda^S (r(l) + \rho_s(l) Q(l))$ reflects the power capacity from swapping stations, including both spare battery stock and in-transit batteries. The third term, $(\lambda^C + \lambda^S)\left(\mathbb{E}[L_z(l)] + \Phi^{-1}(1 - \epsilon^B) \sqrt{\text{Var}(L_z(l))}\right)$, accounts for the expected power loss due to transportation-induced delays and variability. Dividing this total available regulation power by the performance requirement $\eta_z$ yields an upper bound on the regulation capacity $\bar{B}_z(l)$ that can be offered to the grid. This bound will later be integrated into the cost density function in Equation~\eqref{eqn:C_density} as a deduction term, representing the profit from providing frequency regulation services in the centralized charging.

Similarly, we derive the regulation capacity constraint for the benchmark decentralized charging paradigm, where frequency regulation is solely provided by the local spare battery inventory at each swapping station:
\begin{result}
    In the decentralized charging paradigm, the available regulation capacity in the \( z \)th period must satisfy the following condition:
    \begin{align*}
        \bar{B}_z^B(l) \leq \frac{\lambda^S r^B(l)}{\eta_z} = \frac{\lambda^S \left(T^C \bar{\mu}(l) + \Phi^{-1}(1-\epsilon^{BS}) \sqrt{T^C \rho_s(l)} \, \bar{\sigma}(l) \right)}{\eta_z}.
    \end{align*}
\end{result}
This inequality indicates that the maximum regulation capacity $\bar{B}_z^B(l)$ that can be provided in the decentralized paradigm is limited by the total charging and discharging power of local spare batteries and the number of them at swapping stations. The term $\lambda^S r^B(l)$ denotes the maximum power available from decentralized battery reserves, where $r^B(l)$ is the minimum spare inventory level derived in Result~\ref{result:r^B_continuous}. Dividing by the performance requirement $\eta_z$ ensures that a performance threshold $\epsilon^{BS}$ is met under regulation service commitments. The upper bound of this benchmark capacity $\bar{B}_z^B(l)$ will later be integrated into the cost density function in Equation~\eqref{eqn:cost_density_bench} as a deduction term, representing the profit from providing frequency regulation services in the decentralized charging paradigm.

In summary, we have derived the inventory management policies and the corresponding regulation capacities under both decentralized and centralized charging paradigms. Having characterized the operational mechanisms of swapping and charging stations, we proceed to derive the cost density and battery density functions in the following section.

\section{Infrastructure Deployment Modeling and Optimization}
\label{sec:cost_function}

We develop a continuous model to characterize the optimal deployment of battery swapping and charging stations in urban areas. Our objective is to minimize the total amortized cost, which includes infrastructure depreciation, electricity and transportation costs, and accounts for income from frequency regulation services. We analytically derive the cost density and battery density functions for the centralized charging paradigm and the decentralized charging paradigm.

In the centralized charging paradigm, the battery density per unit area around location $l$ is given by  
\begin{align}
\label{eqn:battery_density}
B(l)=R(l)+r(l)+Q(l)\rho_s(l),
\end{align}
where $\rho_s(l)$ denotes the density of swapping stations around location $l$, $R(l)$ is the aggregate required primary stock of all charging stations in a unit area, $r(l)$ is the aggregate re-order point of all swapping stations, and $Q(l)$ is the re-order quantity of each swapping station from the charging station.  
Equation~\eqref{eqn:battery_density} quantifies the total number of batteries required to maintain reliable service in a centralized network.

Let $c^C$ denote the depreciation cost of a charging station, and let $c^R$ denote the depreciation cost of a battery participating in frequency regulation, which is higher than that of a non-regulating battery ($c^R > c^B$). We use $\kappa_n$ to denote the proportion of time in period $z \in \mathcal{Z}$ and $\bar{\mu}_{zn}(l)$ to denote the mean number of EV arrivals at swapping stations per unit area around location $l$ in period $z$ under scenario $n \in \mathcal{N}$.  
Let $c_n^E$ denote the electricity cost per $\text{kW}\cdot\text{h}$ in scenario $n$. Accordingly, the cost density (in $\$ \cdot \mathrm{km}^{-2}$) per unit area around location $l$ is expressed as
\begin{equation}
\begin{array}{rl} \label{eqn:C_density}
C(\rho_c(l),Q(l))=&\sum_{z\in\mathcal{Z}}\sum_{n\in\mathcal N} c_{n}^E \kappa_n \bar{\mu}_{zn}(l) B^C 
+ c^{C}\rho_c(l) 
+ c^{R}\big(R(l)+r(l)+Q(l)\rho_s(l)\big)\\[4pt]
&+\frac{2\sqrt{2}c^T}{3\sqrt{\rho_c(l)}\, Q(l)} 
\sum_{z\in\mathcal Z}\sum_{n\in\mathcal N}\kappa_n \bar{\mu}_{zn}(l)
-\sum_{z\in\mathcal Z} p_z \bar B_z(l),
\end{array}
\end{equation}
where $B^C$ is the battery capacity and $\rho_c(l)$ is the density of charging stations around location $l$. The first term on the right-hand side of Equation~\eqref{eqn:C_density} represents the electricity cost. The second and third terms correspond to the depreciation costs of charging stations and batteries, respectively. The fourth term is the transportation cost, and the last term captures the income from frequency regulation. For any candidate decisions on re-order quantity $Q(l)$ and charging-station density $\rho_c(l)$, a decision maker can use this formulation to quickly evaluate the economic performance of deploying battery-swapping infrastructure around location $l$. Therefore, the unconstrained optimization problem for the battery swapping and charging service operator is 
\begin{equation*}
    \min_{\rho_c(l),Q(l)}\quad C(\rho_c(l), Q(l)).
\end{equation*}
Integrating the cost density over the entire service region $\Omega$ yields the total infrastructural cost $J^{CA}=\int_{l \in \Omega} C(\rho_c(l),Q(l)) \,\mathrm{d}x$. This integral represents the overall objective function in the region. Since the decisions and cost density at each location $l$ are independent of one another, the overall optimization problem can be expressed as
\begin{equation*}
\int_{l \in \Omega}\ \min_{\rho_c(l), Q(l)} C(\rho_c(l),Q(l)) \,\mathrm{d}x, \notag
\end{equation*}
which is equivalent to minimizing the cost at each location $l$ independently. This optimization problem can be efficiently solved because for any location $l$ there are only two decision variables $Q(l)$ and $\rho_c(l)$, and numerical tests indicate that $C(\rho_c(l), Q(l))$ is jointly quasi-convex in $Q(l)$ and $\rho_c(l)$ (see Figure~\ref{fig:sensitivity_deviation}). Besides, the objective function is twice differentiable with either piece of $\phi(Q)$ in \eqref{eq:phi_Q} considered. Therefore, a simple routine with a two-dimensional line search procedure can solve the minimization problem.

For the benchmark decentralized charging paradigm, the battery density around location $l$ is given by
\begin{align}
\label{eqn:battery_density_bench}
B^B(l) &= r^B(l),
\end{align}
and the corresponding cost density is
\begin{align}
\label{eqn:cost_density_bench}
C^{B}(l) &= \sum_{z \in \mathcal{Z}} \sum_{n \in \mathcal{N}} c_{n}^E \kappa_n \bar{\mu}_{zn}(l) B^C + c^I \rho_s(l) + c^R B^B(l) - \sum_{z \in \mathcal{Z}} p_z \bar{B}_z^B(l).
\end{align}
Equation~\eqref{eqn:battery_density_bench} defines the battery density, representing the total number of batteries required to maintain reliable service in a decentralized network per unit area around location $l$.
Equation~\eqref{eqn:cost_density_bench} gives the system cost per unit area around location $l$. Note that this is not an optimization problem, since the stock level is determined by the demand and service level. The first term on the right-hand side of Equation~\eqref{eqn:cost_density_bench} represents the electricity cost for battery charging. The second and third terms correspond to the depreciation costs associated with the construction of swapping stations and the battery inventory, respectively. The last term captures the income reduction from frequency regulation participation, which appears as a negative component. This serves as a benchmark to be compared in the numerical analysis with the centralized charging paradigm. Having derived the continuous cost and battery density functions for both centralized and decentralized charging paradigms, we now examine their practical implications. In the next section, we conduct numerical experiments to analyze how total cost, battery density, and regulation capacity vary with demand scale, charging power, and battery cost.

\section{Numerical Experiments}
\label{sec:numerical_analysis}

In this section, we explore the managerial implications of battery swapping service networks using continuous approximation models. Section~\ref{sec:data_metric} outlines the data sources and the processing of raw data in the frequency regulation market, and introduces four configurations defined by two choices: centralized vs. decentralized charging, and participation vs. non-participation in the frequency regulation market. It also introduces the metrics used to evaluate and compare the performance of these configurations. Section \ref{sec:demand_scale} examines how the configurations perform as demand for battery swapping scales up, while Section \ref{sec:advance_tech} analyzes their performance in response to technological advancements, such as improved charging power and reduced battery costs. Section \ref{sec:flexibility} investigates the sensitivity of cost density to deviations from the optimal re-order quantity and the density of charging stations. Through these numerical experiments and analyses, we derive key managerial insights to support decision-making of stakeholders in battery swapping service.

\subsection{Data and Metrics}
\label{sec:data_metric}

Our location and inventory parameters, including station construction costs, truck transportation, station density, swapping demand, and battery depreciation, are grounded in real-world industry data from \cite{qi2023scaling} and practitioner consultations in Beijing. The amortized costs are $c^S = \$\ 4.46\ /\text{hour}$ for a swapping station, $c^C = \$\ 11.10\ /\text{hour}$ for a charging station, and $c^{I} = \$\ 4.64\ /\text{hour}$ for a swapping station with on-site charging. Battery transportation costs $c^T = \$\ 1.13\ /\text{km}$ with a truck speed $q = 30\ \text{km}/\text{hour}$ and capacity $\bar{Q} = 30$ batteries. Each battery costs $\$\ 7,062.15$, with an eight-year lifespan and a baseline depreciation of $c^B = \$\ 0.10\ /\text{hour}$. Frequent frequency regulation raises the degradation rate from 2.5\% to 7.5\% \citep{sabadini2025does}, shortening lifespan to about three years and increasing depreciation to $c^R = \$\ 0.27\ /\text{hour}$. Charging time per battery is $T^C = 0.78$ hours with a standard fast charger. Electricity prices are $c^E = \$\ 0.223\ /(\text{kW}\cdot\text{h})$ during peak hours (8-10 a.m., 3-8 p.m.) and $\$\ 0.068\ /(\text{kW}\cdot\text{h})$ off-peak \citep{he2022integrated}. Centralized and on-site charging powers are $\lambda^C = 41\ \text{kW}$ and $\lambda^S = 7\ \text{kW}$, respectively. For a single swapping station, hourly swapping demand average $\mu_{z}$ is 5.68 (std. $\sigma_{z}=3.71$) from 0:00 to 8:59 and 14.51 (std. $\sigma_{z}=3.90$) from 9:00 to 23:59. Service levels are set at $\epsilon^S = \epsilon^C = \epsilon^B = 0.03$, and the stockout probability $\epsilon^{BS}$ is calibrated according to Equation~\eqref{eq:epsilon_BS} to ensure comparable service across network designs as proved by \cite{qi2023scaling}:
\begin{equation}
\epsilon^{BS} = \frac{ \mathbb{E} \left[ \left[ \left( N^S(\Delta) - \frac{\rho_c}{\rho_s} R \right)^+ + N^S(T^T) - r \right]^+ \right] }{Q}.
\label{eq:epsilon_BS}
\end{equation}

\begin{table}[h!]
\centering
\caption{Hourly Swap Demands}
\begin{tabular}{ccc}
\toprule
Period & Mean & Standard Deviation \\
\midrule
0:00-8:59  & 5.68  & 3.71 \\
9:00-23:59 & 14.51 & 3.90 \\
\bottomrule
\label{tab:statistic_values}
\end{tabular}
\end{table}

For the frequency regulation component, we use historical data from the PJM website\footnote{https://www.pjm.com/markets-and-operations/ancillary-services}, which includes AGC signals $g_{w,z}$, regulation performance levels $\theta$, and regulation capacity clearing prices $p_z$. For our analysis, we focus on data from August 31, 2024, which is the newest date of data from the PJM website.

To evaluate system performance across key dimensions, we define three core metrics: cost density, battery density, and average frequency regulation capacity, which are associated respectively with different stakeholders such as service providers, environmental stakeholders, and energy providers. Cost density reflects the operational viability and profitability of charging and swapping services. A lower value indicates better economic sustainability. Battery density, although determined by service providers to maximize profitability, also has significant environmental implications. Reducing battery density helps mitigate the adverse environmental impacts associated with battery production, degradation, and recycling processes. Average frequency regulation capacity quantifies the system's contribution to grid stability. A higher frequency regulation capacity signifies a stronger ability to support grid stability, benefiting energy providers. Under a centralized charging paradigm, it can be computed based on Equation~\eqref{eq:B_x}:
\begin{align}
\mathbb{E}\left[B(l)\right] 
&= \frac{\sum_{z\in \mathcal{Z}}\sum_{n\in\mathcal{N}} \kappa_n \mathbb{E}\left[B_{zn}(l)\right]}{|\mathcal{Z}|},
\label{eq:B_x}
\end{align}
where $\mathbb{E}\left[B_{zn}(l)\right] 
= \lambda^C R(l) + \lambda^S \left(r(l) + Q\rho_s(l)\right) - (\lambda^C + \lambda^S) \mathbb{E}\left[L_{zn}(l)\right]$.

In our numerical analysis, we examine two key extensions of battery swapping services: centralized charging and participation in the frequency regulation market. Our baseline is decentralized charging without frequency regulation participation. Accordingly, we conduct a comprehensive comparison across four configurations: (i) decentralized charging, (ii) centralized charging, (iii) decentralized charging with frequency regulation, and (iv) centralized charging with frequency regulation. The results and detailed analysis are presented in the following subsections.

\subsection{Impact of Demand Growth on Charging and Swapping Systems}
\label{sec:demand_scale}

\begin{figure}[htbp]
	\centering
	{\includegraphics[scale=0.45]{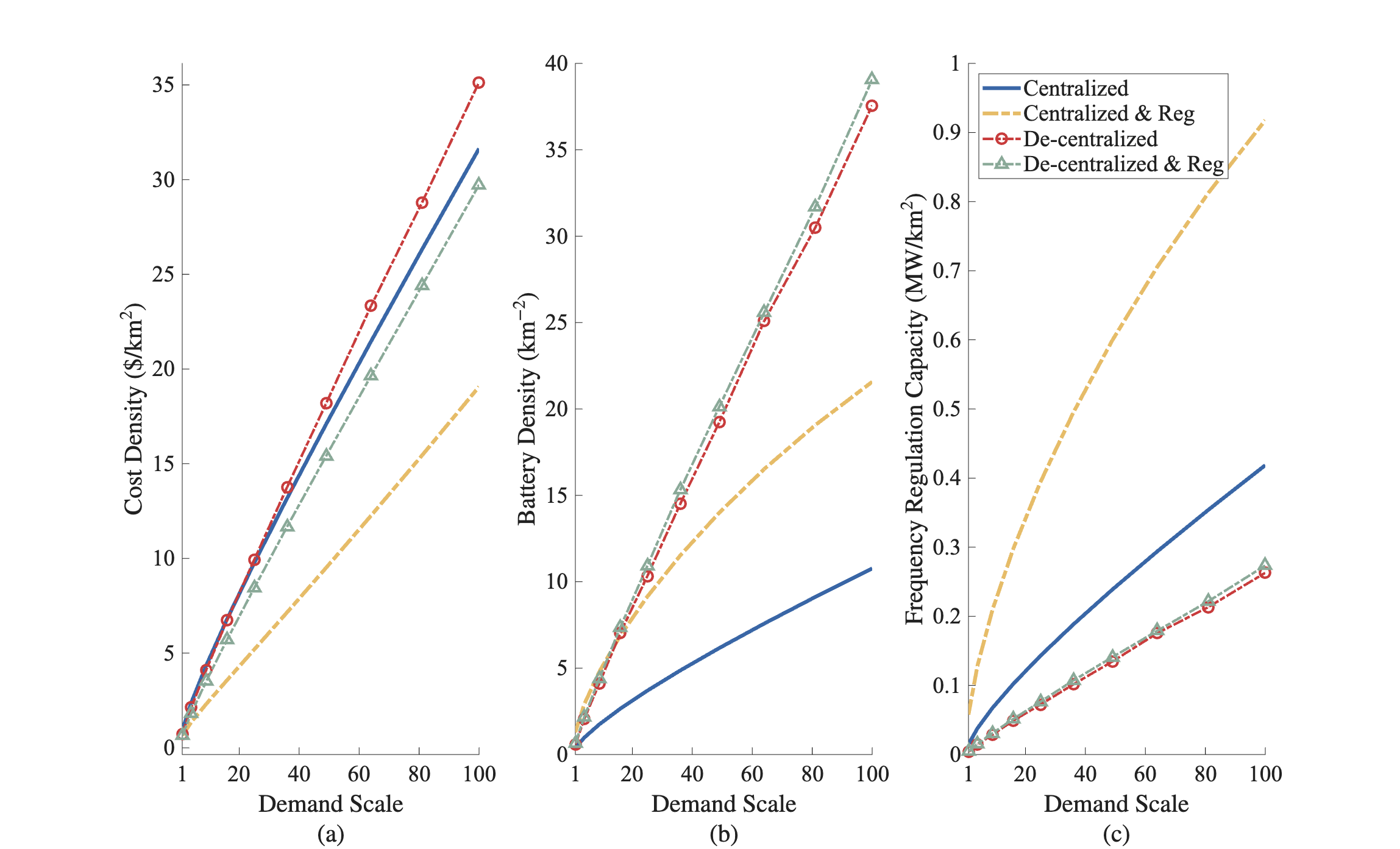}}
	\caption{
    Evolution of (a) the Cost Density; (b) the Battery Density; and (c) the Average Frequency Regulation Capacity as Demand Scales Up Under the Four Configurations}
    \label{fig:cost_battery_density}
\end{figure}

The adoption of EVs has been accelerating rapidly, with over 26 million EVs on the roads by 2022, marking a 60\% increase from 2021 and a fivefold increase since 2018 \citep{KO2024124095}. As a result, the growing number of EVs necessitates the strategic planning and scaling of battery swapping infrastructure, which includes determining the optimal locations and sizes of charging stations to meet the increasing demand \citep{NARASIPURAM2021103225}. Therefore, it is important to determine which charging and swapping service mode can adapt to such scaling up of demand with better performance. In this subsection, we compare the four configurations to find out the efficiency of centralized charging and participation in power frequency regulation market, based on the three metrics. To conduct the comparison, we set centralized charging operates at a higher power level (41 kW), while decentralized charging is constrained to a lower power level (7 kW), which is the same as Level 2 AC charging. For the baseline demand scale, $\rho_{s}$ is set as $0.04~km^{-2}$, which is $10\%$ of the mean density of the gasoline refueling stations in Beijing. Baseline $\bar{\mu}/\rho_{s}$ and $\bar{\sigma}^{2}/\rho_{s}$ are set at the average level observed from the ten swapping stations in our data set, as Table \ref{tab:statistic_values} summarizes. Besides, we set that both the swapping station density $\rho_{s}$ and the demand that each swapping station faces (with mean $\bar{\mu}/\rho_{s}$ and variance $\bar{\sigma}^{2}/\rho_{s}$) grow simultaneously and linearly on a scale from 1 to 10, denoted by $q$. In other words, the demand density $(\bar{\mu}, \bar{\sigma}^{2})$ increases on a scale from 1 to 100. The results and analysis are shown in Figure \ref{fig:cost_battery_density} and below.

\subsubsection{Cost Density Implications}\leavevmode

Cost efficiency is a primary concern for charging station providers, and we evaluate it using cost density, where lower values indicate more cost-effective operations. Figure~\ref{fig:cost_battery_density}(a) plots cost density against demand scale for the four configurations considered in this study: centralized charging and decentralized charging, each with and without participation in the frequency regulation market. Across all configurations, cost density increases sublinearly with demand. Participation in the frequency regulation market shifts each curve downward, and among all cases, centralized charging with frequency regulation consistently achieves the lowest cost density over the entire demand range.

This pattern can be traced to two main mechanisms. First, joining the frequency regulation market directly reduces net cost density under both centralized and decentralized configurations, because it generates an additional revenue stream that offsets the costs. The magnitude of this reduction is greater for centralized charging. This is because higher charging power at centralized facilities allows more regulation capacity to be bid per unit of battery, and pooled batteries enable the system to offer a greater regulation capacity. Second, decentralized charging yields higher cost density because its lower charging power at each site necessitates larger battery reserves to avoid stockouts, which raises inventory-related costs. This is because lower charging power lengthens effective replenishment times, so more batteries are required to maintain the same service level. In addition, low charging power reduces the regulation capacity that each battery can provide, thus reducing the income from participation in frequency regulation.

A related scaling effect explains the observed sublinear growth of cost density. As the demand variance $\frac{\bar \sigma^2}{\rho_s}$ scale up proportionally, uncertainty is aggregated and pooled across stations, which means the required safety stock grows more slowly than demand variance itself. This reduction in relative variability leads to an economy of scale.

From a managerial perspective, these results imply that service providers seeking to minimize costs should prioritize building centralized charging stations to serve decentralized swapping stations. Furthermore, regardless of the spatial configuration chosen, service providers should seek to participate in the frequency regulation market. The resulting income can meaningfully offset the substantial investment and operating costs associated with charging and swapping infrastructure.

\subsubsection{Battery Density Implications} \leavevmode

Figure~\ref{fig:cost_battery_density}(b) plots battery density against demand scale for the four configurations considered in this study. It shows that participation in the frequency regulation market increases battery density under both centralized and decentralized charging configurations. At the same time, for any given demand level, decentralized charging requires a substantially larger battery stock than centralized charging. These two effects are persistent across the entire demand range, and the gap between decentralized and centralized configurations widens as demand scales up.

This pattern can be traced to two main mechanisms, similar to those in the cost density analysis.
First, participation in the frequency regulation market allows service providers to earn extra revenue. They can bid more capacity into the market when they hold more batteries. As a result, they are incentivized to keep a larger stock than what is required to meet swapping demand alone. Second, decentralized charging increases effective demand variability at each station. The lower charging power prolongs replenishment times and reduces pooling of fluctuations across the network. This greater variability forces operators to maintain higher safety stock to avoid stockouts. It leads to consistently higher battery density than centralized charging, regardless of regulation participation.

These results imply that while frequency regulation participation offers financial benefits, it also introduces additional inventory. To minimize battery-related environmental and economic burdens, centralized charging remains the preferable strategy, as it achieves significant reductions in battery stock while still enabling profitable market participation. The combination of centralized charging and frequency regulation participation continues to exhibit a favorable trade-off balance between economic profit and environmental performance.

\subsubsection{Average Frequency Regulation Capacity Implications} \leavevmode

Figure~\ref{fig:cost_battery_density}(c) plots the average frequency regulation capacity, which measures the system's ability to contribute to grid stability, for the four configurations considered in this study. Across all configurations, regulation capacity increases with demand. Centralized charging with frequency regulation consistently achieves the fastest growth, followed by centralized charging without frequency regulation. In contrast, decentralized charging, regardless of frequency regulation participation, exhibits slower growth, underscoring the superior regulation potential of centralized charging strategies.

Both centralized charging and participation in the frequency regulation market contribute to higher regulation capacity, with their combination yielding the greatest gains. Centralized charging enables substantially higher charging power ($\lambda^C$), allowing more regulation capacity to be bid into the market per unit of battery. Meanwhile, participation in the frequency regulation market incentivizes service providers to increase battery stock in order to provide greater regulation capacity and thus generate additional revenue. The interaction of these two mechanisms explains why both centralized charging and frequency regulation independently increase regulation capacity, and why their combination yields the highest capacity across the entire demand range.

In summary, investing in centralized charging infrastructure and actively participating in the frequency regulation market can significantly enhance an service provider's ability to stabilize grid frequency. This strategy strengthens the system's role as a flexible grid resource while also generating ancillary service revenues, benefiting both the service provider and the power system.

\subsubsection{Summary of Comparison of Four Configurations with Demand Scaling up} \leavevmode

\begin{figure}[htbp]
	\centering
	{\includegraphics[scale=0.35]{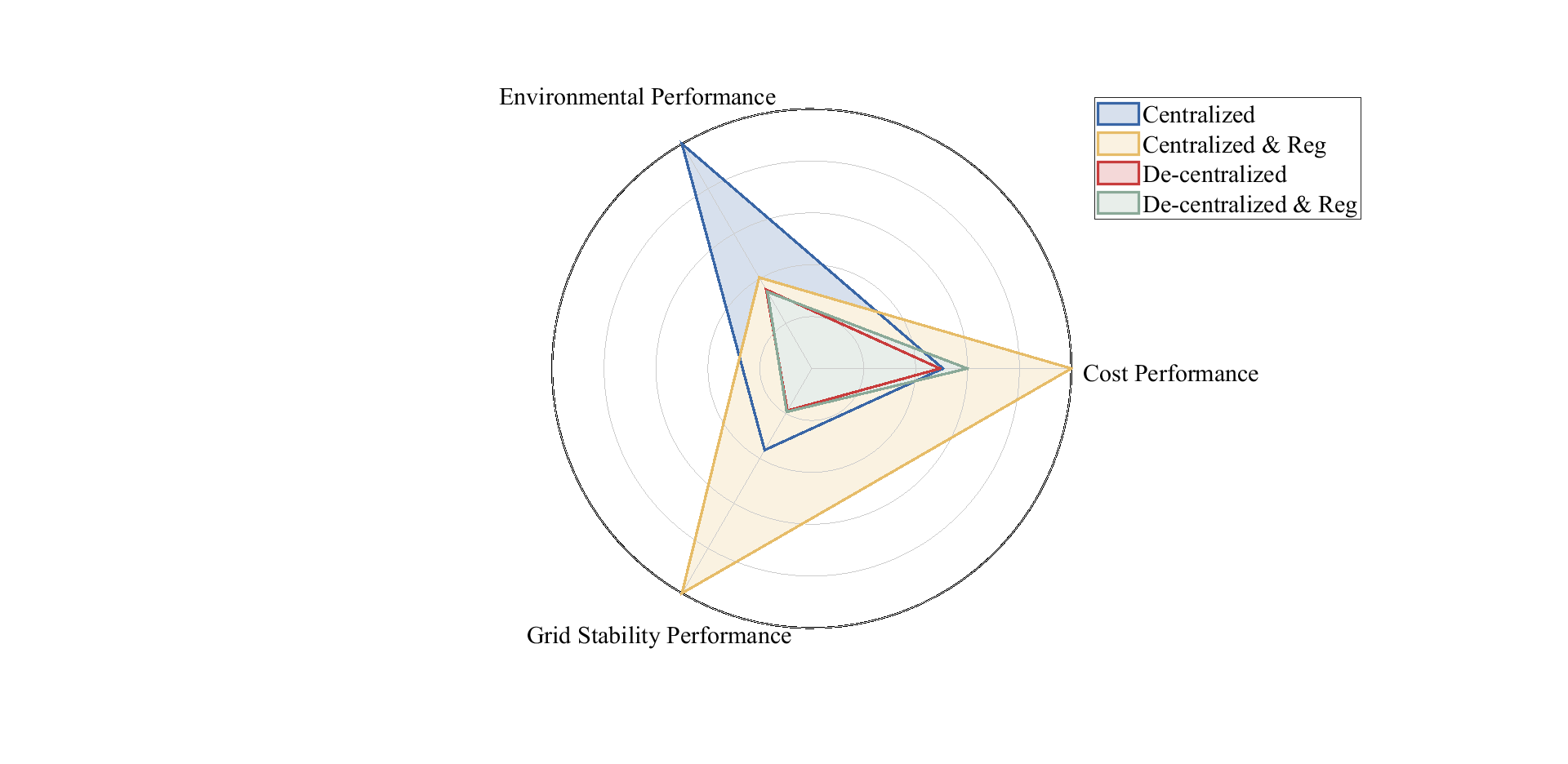}}
	\caption{
    Comparative Performance Evaluation of Four Configurations on Three Key Dimensions at $q=5$ Using Radar Chart}
    \label{fig:radar_demand}
\end{figure}

Figure \ref{fig:radar_demand} shows a radar chart comparing the performance of four configurations across three dimensions at scale $q=5$. The dimensions include Cost Performance and Environmental Performance, measured inversely by cost density and battery density respectively, along with Grid Stability Performance. All three metrics are normalized to a [0,1] scale where values closer to 1 indicate better performance. The results clearly demonstrate that the configuration combining centralized charging with frequency regulation delivers the strongest overall performance, followed by centralized charging alone. Both decentralized charging configurations show relatively poor performance, though the version incorporating frequency regulation performs marginally better than its counterpart without regulation.

In summary, we list our key findings of this subsection in which we compare the different performance of the four configurations under different demand scales:

\begin{quote}
\textbf{Finding 1.}
\begin{enumerate}[label=(\roman*)]
    \item \textbf{Centralized Charging} excels across all metrics with demand scaling up. With higher charging power (41 kW), centralized charging lowers cost density and battery density, improving cost efficiency and environmental performance by minimizing the need for large battery reserves. It also significantly enhances frequency regulation capacity, improving grid stability as demand grows.
    \item \textbf{Participation in Frequency Regulation} results in lower cost density and higher frequency regulation capacity but higher battery density, which means it improves cost efficiency and grid stability but slightly worsen environmental performance with demand scaling up. Nevertheless, the worse environmental impact is relatively minor compared to the overall environmental benefits of centralized charging.
    \item \textbf{Combination both approaches} results in lower cost density and higher frequency regulation capacity, which means it improves performance on cost efficiency and grid stability compared to standalone implementations. Besides, its battery density representing environmental performance is between the the standalone implementations of the two approaches.
\end{enumerate}
\end{quote}

\subsection{Impact of Technological Advancements on Charging and Swapping Systems}
\label{sec:advance_tech}

Technological advancements have significantly improved the efficiency of EV charging and swapping services, particularly through innovations in charging power and battery technology. As charging power increases, the required charging time decreases, reducing station occupancy and lowering the need for large battery reserves. This optimization helps minimize overall infrastructure costs while enhancing service efficiency \citep{Ahn2015}. Moreover, higher charging power can contribute to better grid stability by enabling more flexible integration of renewable energy sources. 

Technological advancements also lower the cost of battery production and maintenance. Studies project substantial reductions in battery prices over time. For example, battery prices are expected to decrease by 10\% in 5 years, 28\% in 10 years, 48\% in 15 years, and 54\% in 20 years \citep{lemme2019optimization}. With costs of battery decreasing, it will help improve the cost efficiency of battery swapping service providers, while also changing their decisions about the building of stations and stocking of batteries.

In this section, we explore whether centralized charging and participation in the frequency regulation market can better leverage higher charging power and lower battery cost. Specifically, we analyze the impact on EV charging and swapping systems from three perspectives: cost efficiency, environmental impact, and frequency regulation capacity. This analysis aims to assess the role of centralized charging and frequency regulation in optimizing the benefits of advanced technologies. Analysis in this part are based on a moderate future demand level ($q=5$).

\subsubsection{Cost Density Implications} \leavevmode

\begin{figure}[htbp]
	\centering
	{\includegraphics[scale=0.45]{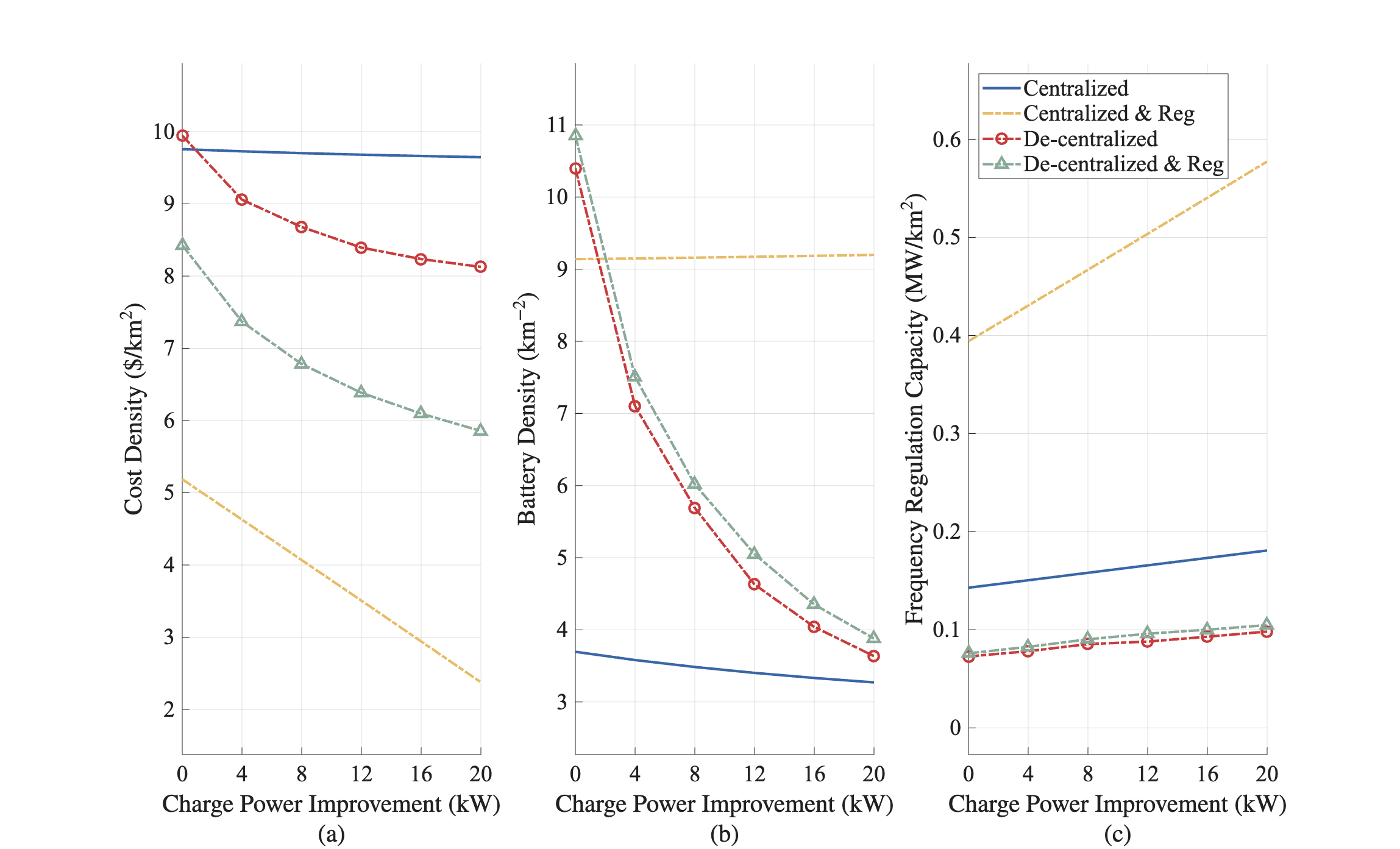}}
    \caption{
    Evolution of (a) the Cost Density; (b) the Battery Density; and (c) the Average Frequency Regulation Capacity as Charging Power Increasing from Baseline Power Level (41 kW, Centralized Charging and 7 kW, Decentralized Charging) Under the Four Configurations  \label{fig:charge_power}}
\end{figure}

\paragraph{\textbf{Increasing charging power:}} Higher charging power generally leads to a decrease in cost density across all configurations, indicating improved cost efficiency as power capacity increases, as illustrated in Figure~\ref{fig:charge_power}(a). Among the four configurations considered, the ranking of cost efficiency from highest to lowest under increasing charging power is: centralized charging with frequency regulation, decentralized charging with frequency regulation, decentralized charging without regulation, and centralized charging without regulation.

In the absence of frequency regulation, the marginal cost savings of increasing charging power diminish for centralized charging, making it less favorable than decentralized charging as power continues to rise. This occurs because the primary advantage of centralized charging - faster battery charging leading to shorter charging time ($T^C$) - gradually weakens at high power levels. When charging power is sufficiently high, decentralized charging can also recharge batteries quickly. At this point, the transportation time ($T^T$) inherent to centralized charging becomes the dominant factor and effectively sets a lower bound on battery replenishment time, limiting further cost reductions from increasing charging power. Decentralized charging, not being constrained by this transportation delay, thus gains relatively more benefit from increased charging power when frequency regulation is not considered.

However, participation in the frequency regulation market plays a crucial role in reversing this trend, making centralized charging more cost-effective than decentralized charging as charging power increases. This pattern emerges because centralized charging combined with frequency regulation benefits most from higher charging power, which enables greater regulation capacity and consequently higher profitability. By contrast, decentralized charging configurations are limited by lower station charging power and less efficient aggregation of regulation capacity.

These results suggest that as grid power capacity improves, the limiting effect of transportation time on centralized charging reduces the benefits of increasing charging power alone. However, participation in the frequency regulation market introduces additional revenue incentives that restore the advantage of centralized charging over decentralized alternatives. Therefore, battery swapping service providers can achieve the greatest gains by combining centralized charging infrastructure with frequency regulation market participation.

\begin{figure}[htbp]
	\centering
	{\includegraphics[scale=0.45]{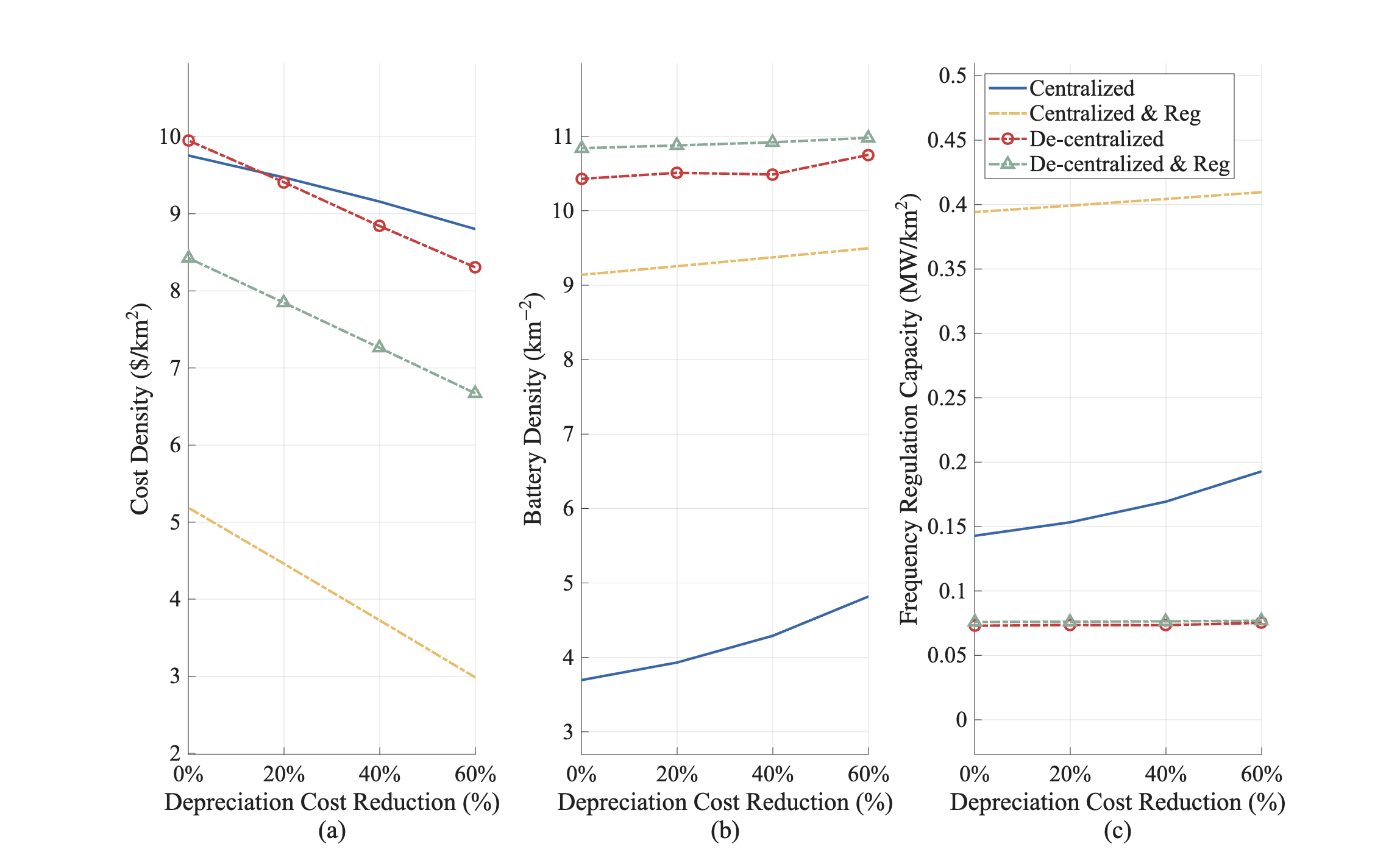}}
    \caption{
    Evolution of (a) the Cost Density; (b) the Battery Density; and (c) the Average Frequency Regulation Capacity as Battery Depreciation Cost Decreasing from Baseline Under the Four Configurations  \label{fig:depreciation_cost}}
\end{figure}

\paragraph{\textbf{Reducing battery costs:}} Figure~\ref{fig:depreciation_cost}(a) shows that a reduction in battery depreciation costs leads to a decrease in overall cost density. However, when participation in the frequency regulation market is not considered, centralized charging no longer demonstrates a clear cost advantage over the baseline decentralized charging approach. This indicates that the economic benefit of centralized charging diminishes as battery costs decline.

This shift can be explained by two main factors. First, the primary advantage of centralized charging lies in its ability to reduce battery inventory requirements compared to decentralized charging. This advantage becomes less significant when batteries are cheaper. Second, centralized charging incurs substantial infrastructure costs for dedicated charging stations, which are not offset by savings from lower battery costs. Consequently, without additional revenue streams, the economic viability of centralized charging weakens as battery costs fall compared to decentralized charging.

However, participation in the frequency regulation market mitigates this challenge by providing a steady source of revenue from grid services, which is largely unaffected by battery cost reductions. Moreover, combining centralized charging with frequency regulation enables the system to maintain smaller battery inventories while simultaneously earning higher profits through greater charging power and increased regulation capacity.

From a managerial perspective, these results suggest that as battery costs continue to decrease, service providers should actively seek to diversify their revenue streams, particularly by participating in frequency regulation markets, to sustain the economic benefits of centralized charging infrastructure.

\subsubsection{Battery Density Implications} \leavevmode

\paragraph{\textbf{Increasing charging power:}}  
Figure~\ref{fig:charge_power}(b) illustrates that under the decentralized charging paradigm, battery density decreases as charging power increases, leading to a reduction in environmental impact. In centralized charging without participation in the frequency regulation market, battery density remains low and stable across different power levels. Conversely, centralized charging combined with frequency regulation exhibits a consistently high and stable battery density, reflecting the economic incentive to maintain larger battery inventories for capacity regulation purposes. Both decentralized charging with or without frequency regulation, the battery density decreases fast as charging power increases.

Participation in the frequency regulation market encourages battery swapping service providers to hold larger battery stocks to enhance regulation capacity and increase revenue, which results in higher battery density compared to setups without regulation. Meanwhile, centralized charging itself promotes environmental benefits by improving battery management efficiency and reducing inventory needs relative to decentralized charging. However, when combined with frequency regulation, these environmental gains are partially offset due to the increased battery stocking required.

For the decentralized charging paradigm, battery density decreases rapidly because shorter charging times reduce the need to maintain large battery inventories to prevent stockouts. Participation in the frequency regulation market also leads to an increase in battery density under decentralized charging compared to the case without frequency regulation; however, this increase is relatively modest. This is because service providers gain limited benefits from decentralized charging in the frequency regulation market and thus have less incentive to significantly raise battery inventory levels.

As charging power increases, the battery density under the decentralized charging paradigm gradually approaches that of centralized charging. However, participation in the frequency regulation market consistently leads to relatively higher battery density, especially in the centralized charging configuration. This creates a trade-off with cost density. Therefore, service providers should carefully balance these considerations when selecting the system configuration as charging power increases.

\paragraph{\textbf{Reducing battery costs:}}  
Figure \ref{fig:depreciation_cost}(b) shows that as battery costs decrease, the battery density remains stable in the baseline decentralized charging paradigm. When participating in the frequency regulation market, the battery density becomes slightly higher than the baseline due to service providers stocking more batteries to increase profits. In contrast, the centralized charging mode reduces battery inventory levels while improving environmental performance through higher charging efficiency. The combined approach shows intermediate environmental benefits - better than using only frequency regulation but not as good as using only centralized charging. Besides, the battery density of centralized charging paradigm increases faster compared to the stable battery density of decentralized charging.

As battery costs decrease, the battery density and their relative rankings across the four configurations remain largely unchanged. This suggests that battery cost has a limited impact on overall inventory strategies. Instead, inventory decisions are more decisively influenced by market factors such as meeting customer demand and participation in the frequency regulation market, which play a more critical role in shaping battery inventory policies.

\subsubsection{Average Frequency Regulation Capacity Implications} \leavevmode
\paragraph{\textbf{Increasing charging power:}}
As illustrated in Figure \ref{fig:charge_power}(c), frequency regulation capacity remains relatively stable in decentralized charging configurations as charging power increases. Centralized charging itself can improve the capacity, but still remains relatively stable with increasing charging power. However, when centralized charging participates in frequency regulation, its regulation capacity increases significantly as charging power rises. This reinforces the advantage of centralized charging in providing greater grid stability. Therefore, we can summarize that both centralized charging mode and participation in frequency regulation market can help improve frequency regulation capacity and benefit the power company, and it will be the most beneficial to adopt the two methods together.

\paragraph{\textbf{Reducing battery costs:}} 
Figure \ref{fig:depreciation_cost}(c) reveals significant differences in regulation capacity between charging architectures. The decentralized configurations demonstrate limited regulation capacity, including the case with frequency regulation participation, primarily constrained by their lower charging power. In contrast, centralized charging achieves substantially higher performance, with its regulation capacity further enhanced through frequency market participation. These results clearly highlight the operational advantages of centralized architectures in providing grid regulation services.

\subsubsection{Summary of Comparison of Four Configurations With Advanced Technologies} \leavevmode

With charging power improving and battery cost decreasing, the two approaches, centralized charging mode and participation in frequency regulation market performs significantly on the three metrics. We summarize the performance of them on different metrics with charging power increasing and battery cost decreasing in Table \ref{tab:advancements_compare}.

\begin{table}[htbp]
    \small
    \centering
    \caption{Performance of Centralized Charging and Frequency Regulation Participation Under Technological Advancements Across Different Charging Power and Battery Cost (q=5)}
    \renewcommand{\arraystretch}{1.5}
    \setlength{\tabcolsep}{10pt}
    \begin{tabular}{m{2.5cm}m{2cm}m{2.7cm}m{2.7cm}m{3.2cm}}
    \toprule
    \textbf{Advancements} & \textbf{Approach} & \textbf{Cost Density} & \textbf{Battery Density} & \textbf{Frequency Regulation Capacity} \\
    \midrule
    \multirowcell{3}[0ex]{\makecell[c]{Charging Power \\ Increasing}} 
    & Centralized Charging & $[-1.89\%, +18.96\%]$ & $[-64.46\%, -10.03\%]$ & $[+84.18\%, +96.23\%]$ \\
    \cline{2-5}
    & Participation in Frequency Regulation & $[-27.95\%, -15.25\%]$ & $[+4.36\%, +9.10\%]$ & $[+4.36\%, +9.10\%]$ \\
    \cline{2-5}
    & Combined Approach & $[-70.77\%, -47.84\%]$ & $[-12.11\%, +153.00\%]$ & $[+441.48\%, +488.16\%]$ \\
    \midrule
    \multirowcell{3}[0ex]{\makecell[c]{Battery Cost \\ Decreasing}} 
    & Centralized Charging & $[-1.96\%, +5.98\%]$ & $[-64.56\%, -55.18\%]$ & $[+95.69\%, +156.30\%]$ \\
    \cline{2-5}
    & Participation in Frequency Regulation & $[-19.69\%, -15.33\%]$ & $[+2.12\%, +4.12\%]$ & $[+2.12\%, +4.12\%]$ \\
    \cline{2-5}
    & Combined Approach & $[-64.08\%, -47.88\%]$ & $[-12.36\%, -10.61\%]$ & $[+439.98\%, +450.84\%]$ \\
    \bottomrule
    \end{tabular}
    \label{tab:advancements_compare}
\end{table}

In summary, we list our key findings of this subsection in which we compare the different performance of the four configurations under different charging power:

\begin{quote}
\textbf{Finding 2.}

Facing technological advancements including higher charging power and lower battery cost,
\begin{enumerate}[label=(\roman*)]
    \item \textbf{Centralized Charging} results in lower battery density and higher frequency regulation capacity, which means better performance on environment and grid stability. However, it is not preferable for battery swapping service providers to adopt centralized charging mode only, as this approach performs worse on cost efficiency.
    \item \textbf{Participation in Frequency Regulation} results in lower cost density and higher frequency regulation capacity, which means better cost efficiency and greater contribution to grid stability, but higher battery density which means worse environmental performance.
    \item \textbf{Combining both approaches} significantly reduces cost density and enhances frequency regulation capacity compared to implementing either approach alone, though at the expense of higher battery density representing reduced environmental performance compared to centralize charging alone.
\end{enumerate}
\end{quote}

\subsection{Flexibility of re-order Quantity and Charging Station Density}
\label{sec:flexibility}

\begin{figure}[htbp]
    \centering
    \begin{subfigure}[t]{0.48\textwidth}
        \centering
        \includegraphics[width=\linewidth]{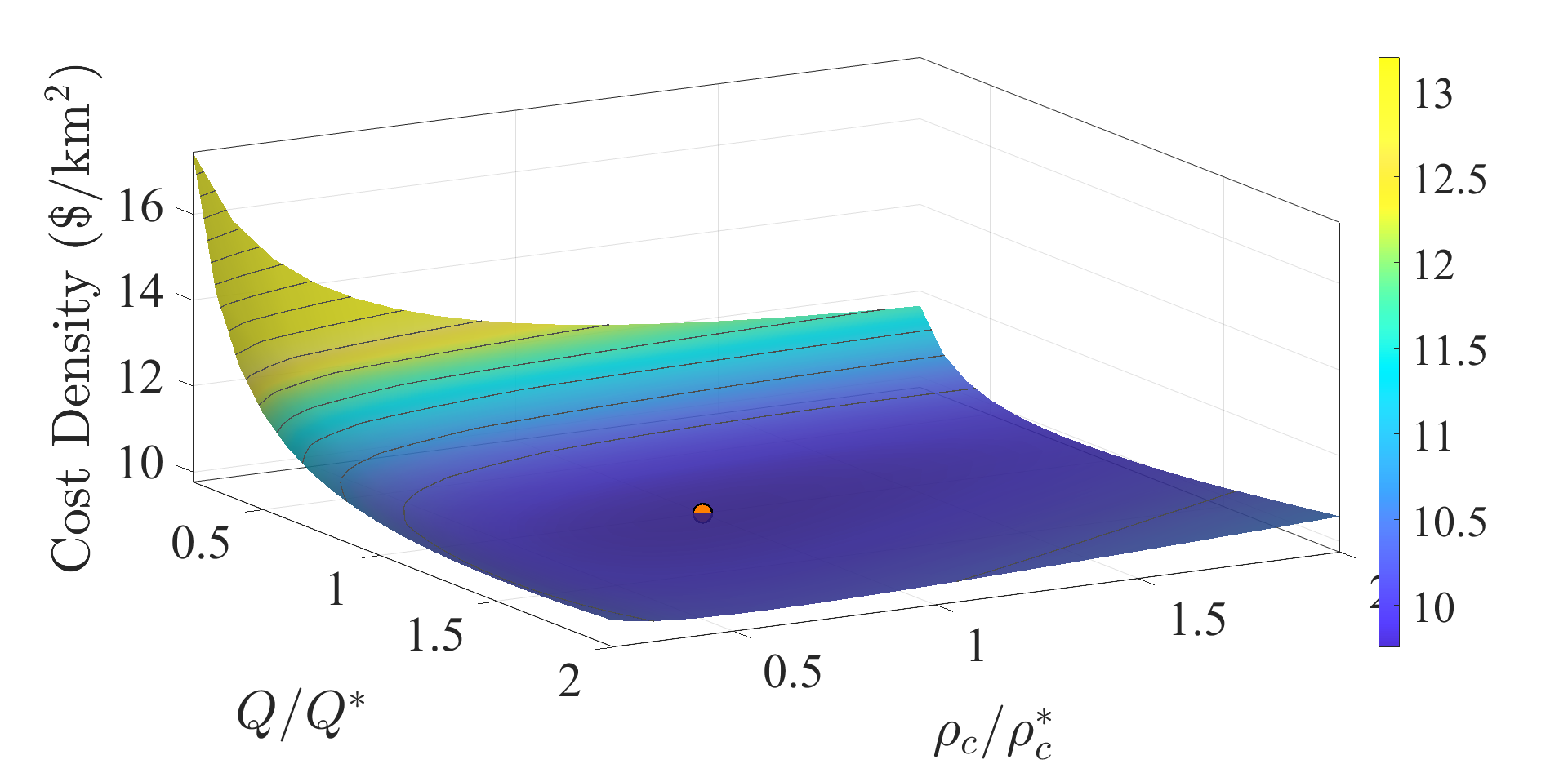}
        \caption{}
        \label{fig:decision_d1}
    \end{subfigure}
    \hfill
    \begin{subfigure}[t]{0.48\textwidth}
        \centering
        \includegraphics[width=\linewidth]{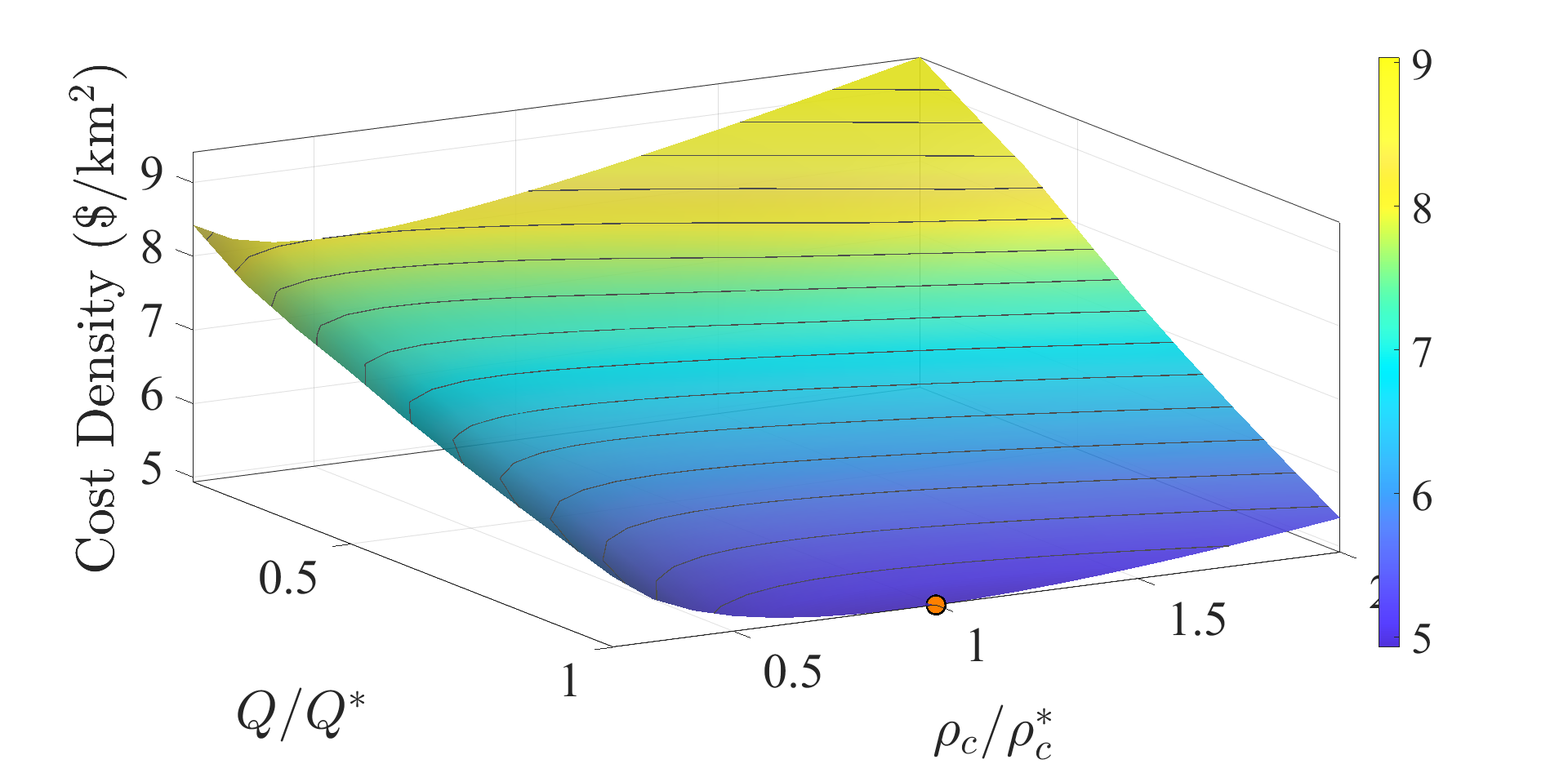}
        \caption{}
        \label{fig:decision_d2}
    \end{subfigure}
    \caption{
    Sensitivity of the Cost Density to the Deviation from
    the Optimal Solution. The densities of swapping stations and demand are set as five times the baseline
    values. The orange dots are optimal solutions.  \label{fig:sensitivity_deviation}}
\end{figure}

In this part, we discuss the flexibility of the two decisions in centralized charging modes: re-order quantity and charging station. Will deviating from the optimal re-order quantity significantly increase the total cost? This question is motivated by the reality where the optimal battery re-order quantity may not necessarily equal to one truck load. The trucking cost may also be tiered so that only several discrete values of $Q$ are implementable. Meanwhile, it is environmentally desirable to reduce $Q$ to reduce the battery stock level. It turns out that the cost density is a remarkably flat function in $Q$ around its optimal value in the centralized charging mode without frequency regulation, as shown in Figure \ref{fig:sensitivity_deviation}(a). Therefore, in this case re-order quantity can be relatively less than the optimal solution, which serves as a bargain zone between agencies of battery swapping service providers and agencies of environmental departments. However, when the participation in frequency regulation is considered, battery swapping service providers tend to maximize re-order quantity to the upper bound of stocking spaces of stations or trucks, as they benefit from this, and also a little decrease in optimal re-order quantity will cause a lot of loss for service providers, as the market profit earned from regulation market per battery stocked is very high.

In this section, we examine the flexibility of two key decisions in centralized charging modes: re-order quantity and charging station selection. A crucial question arises - how sensitive is the total cost to deviations from the optimal re-order quantity? This question is particularly relevant in practice, as the optimal battery re-order quantity may not necessarily align with a full truckload. Moreover, trucking costs are often tiered, meaning only certain discrete values of $Q$ are feasible. At the same time, reducing $Q$ is environmentally preferable as it helps lower battery stock levels.

Interestingly, in the centralized charging mode without frequency regulation, the cost density function remains remarkably flat around its optimal re-order quantity, as illustrated in Figure \ref{fig:sensitivity_deviation}(a). This suggests that slight reductions in re-order quantity do not significantly impact total costs, providing a ``bargaining zone" between battery swapping service providers and environmental agencies. However, the situation changes when participation in frequency regulation is considered. In this case, battery swapping service providers tend to maximize their re-order quantity, pushing it to the upper stocking limits of stations or trucks. The reason is that the financial gains from the frequency regulation market per stocked battery are substantial, making even a slight reduction in the optimal re-order quantity result in significant losses. Thus, while flexibility exists in the absence of frequency regulation, the economic benefits of regulation drive service providers toward maximizing re-order quantities whenever possible.

Urban space presents inherent complexities, with regulations and permanent structures often restricting battery swapping service providers from selecting the ideal locations for deploying charging stations. However, a fortunate observation from Figures \ref{fig:sensitivity_deviation}(a) and (b) is that the cost density function remains remarkably flat with respect to charging station density. This implies that service providers have considerable flexibility in adjusting the placement of charging stations to better navigate urban constraints - factors that are beyond the scope of this paper - without incurring significant additional costs.

\begin{quote}
\textbf{Finding 3.}

In the centralized charging mode, the density of charging stations exhibits high flexibility, accommodating the complexities and spatial constraints of urban space. Additionally, the flexibility in the optimal re-order quantity provides room for negotiation between service providers and environmental stakeholders. However, when participating in the frequency regulation market, cost density becomes highly sensitive to deviations from the optimal re-order quantity. This sensitivity may pose challenges for service providers seeking to mitigate environmental impact by reducing re-order quantities.
\end{quote}

\section{Conclusion}
\label{sec:conclusion}

With the rapid development of electric vehicles and their charging and swapping services, this paper is motivated by the challenges faced by urban EV battery swapping systems in balancing customer demand, infrastructure cost, environmental sustainability, and grid compatibility. To understand and address these challenges, we examine four operational configurations that combine decentralized versus centralized charging with participation in frequency regulation markets. We develop an integrated modeling framework that jointly captures location planning, battery inventory, and grid interaction, and employ a continuous approximation approach to enable tractable optimization of high-dimensional location-inventory-grid coupling problem while preserving key managerial insights, which is a task prohibitively complex for discrete models. Calibrated on real-world data, our numerical experiments reveal three findings: (1) Centralized charging markedly reduces battery stock requirements and cost density compared with decentralized charging, especially when charging power is high; (2) Participation in frequency regulation provides significant cost offsets and improves grid stability, though it modestly increases battery inventory; and (3) Combining centralized charging with frequency regulation achieves the most balanced performance across cost, environmental, and grid-stability dimensions.

This paper offers an analytical foundation for understanding how infrastructure configuration and grid integration jointly shape the economic of large-scale battery swapping networks. It also demonstrates the value of continuous approximation methods for exploring complex, coupled infrastructure problems in a scalable manner. Many research opportunities remain for extending this work, such as incorporating stochastic travel times in battery transport, modeling heterogeneous EV battery types and states of charge, and considering dynamic pricing in both swapping services and frequency regulation markets. Moreover, integrating renewable energy supply variability, local microgrid interactions, and battery-to-grid discharge options could yield deeper insights into the mobility-energy nexus, supporting the evolution of cost-efficient, environmentally sustainable, and grid-friendly urban EV infrastructures.

\section*{CRediT authorship contribution statement}
 
Wenqing Ai: Conceptualization, Funding acquisition, Methodology, Writing - original draft, Writing - review \& editing. 
Hanyu Cheng:  Writing - original draft, Writing - review \& editing, Methodology, Validation, Visualization.
Wei Qi: Conceptualization, Supervision, Resources, Funding acquisition, Writing - review \& editing.

\section*{Acknowledgement}

The authors gratefully acknowledge the specific research fund of the National Natural Science Foundation of China and China Postdoctoral Science Foundation. Wenqing Ai would like to acknowledge the support from the National Natural Science Foundation of China [Grants 72301027] and China Postdoctoral Science Foundation [Grants 2023M730215, 2025T180213]. Wei Qi acknowledges the support from the National Natural Science Foundation of China [Grants 72242106, 72521001].

\def\bibfont{\fontsize{11}{14}\selectfont}

\normalem
\bibliographystyle{informs2014}
\bibliography{ref}

\end{document}